\documentclass[final,1p,times]{elsarticle}
\usepackage{amsmath,mathrsfs,amssymb,amsthm,dsfont,bm,amsfonts,booktabs}        
\usepackage{color,graphicx,subfigure,epsfig}
\allowdisplaybreaks
\newtheorem{thm}{Theorem}
\newtheorem{lem}{Lemma}
\newtheorem{assumption}{Assumption}
\usepackage{multirow}
\newtheorem{alg}{Algorithm}
\newcommand{\p}{\partial}
\newcommand{\f}{\frac}

\newcommand{\up}{\textup}

\newcommand{\ds}{\displaystyle}

\journal{*******************}
\begin{document}
\begin{frontmatter}
\title{Numerical simulation of wormhole propagation with the mixed hybridized discontinuous Galerkin finite element method}

\author[a]{Jiansong Zhang\corref{zhang}}
\ead{jzhang@upc.edu.cn}
\author[b]{Jiang Zhu}
\ead{jiang@lncc.br}

\author[a]{Yiming Wang}
\ead{1164324059@qq.com}
\author[a]{Yanyu Liu}
\ead{1872371250@qq.com}

\author[a]{Hui Guo}
\ead{sdugh@163.com}

\address[a]{College of Science,
        China University of Petroleum, Qingdao 266580, China}
\address[b]{Laborat\'orio Nacional de Computa\c{c}\~ao Cient\'{\i}fica, MCTI 
        Avenida Get\'ulio Vargas 333, 25651-075 Petr\'opolis, RJ, Brazil}

\cortext[zhang]{Corresponding author. Jiansong Zhang}

\begin{abstract}
The acid treatment of carbonate reservoirs is a widely employed technique for enhancing the productivity of oil and gas reservoirs. In this paper, we present a novel combined hybridized mixed discontinuous Galerkin (HMDG) finite element method to simulate the dissolution process near the wellbore, commonly referred to as the wormhole phenomenon. The primary contribution of this work lies in the application of hybridization techniques to both the pressure and concentration equations. Additionally, an upwind scheme is utilized to address convection-dominant scenarios, and a ``cut-off" operator is introduced to maintain the boundedness of porosity. Compared to traditional discontinuous Galerkin methods, the proposed approach results in a global system with fewer unknowns and sparser stencils, thereby significantly reducing computational costs. We analyze the existence and uniqueness of the new combined method and derive optimal error estimates using the developed technique. Numerical examples are provided to validate the theoretical analysis.
\end{abstract}

\end{frontmatter}


\section{Introduction}
\setcounter{equation}{0}
In this article, we will consider a new combined method to simulate the  wormhole propagation which is governed by the following  partial differential equations (see \cite{rui0,zhang1,zhang2023}): 
\begin{equation} \label{e1}
\begin{split}
&(\up{a})\quad\f{\p \phi}{\p t}+\nabla\cdot\mathbf{u}=f,\quad\mathbf{u}=\ds-\f{\kappa(\phi)}{\mu}\nabla p,\quad x\in\Omega,\ t \in(0,T],\\
&(\up{b})\quad\f{\p }{\p t}(\phi c_f)+\nabla\cdot(\mathbf{u}c_f-\phi D\nabla c_f)=\kappa_ca_v(c_s-c_f)+f_Pc_f+f_Ic_I, \quad x\in\Omega,\ t \in(0,T],\\
&(\up{c})\quad\f{\p \phi}{\p t}=\ds\f{\alpha \kappa_c a_v (c_f-c_s)}{\rho_s},\quad x\in\Omega,\ t \in(0,T],
\end{split}
\end{equation}
and the  initial-boundary conditions
\[
\left\{\begin{aligned}
&\phi(x,0)=\phi_{0}(x),\quad c_{f}(x,0)=c^{0}_{f}(x),\quad x\in \Omega,\\
&\mathbf{u}\cdot \mathbf{n}=0,\quad \phi \mathbf{D} (\mathbf{u})\nabla c_f \cdot \mathbf{n}=0,\quad x\in \partial\Omega,\quad 0\le t \le T,
\end{aligned}\right.
\]
where $\Omega$ is a polygonal bounded domain in $\mathsf{R}^d(d=2,\ 3)$;  $p$ denotes the pressure, and $\mathbf{u}$ is the Darcy velocity; $\mu$ is the fluid viscosity, and the source function $f =  f_P+f_ I $, here $f_P$ and
$f_I$ are production and injection rates, respectively;   
The diffusion coefficient  $D(\mathbf{u})=d_m I+|\mathbf{u}|(d_l\up{E}(\mathbf{u})+d_t\up{E}^{\perp}(\mathbf{u}))$,   $\up{E}(\mathbf{u})= [u_{i}u_j/|\mathbf{u}|^2]$ and $\up{E}^{\perp}(\mathbf{u})=I-\up{E}(\mathbf{u})$; $d_{m}$, $d_{l}$, and $d_{t}$ denote the molecular diffusion and dispersion coefficients, respectively;  $a_v$ is the interfacial area available for reaction per unit volume of the medium; $\kappa_c$ is the local mass-transfer coefficient; $c_f$ is the cup-mixing concentration of the acid in the fluid phase, $c_I$ is the injected concentration, $c_s$ is the concentration of the acid at the fluid-solid interface, and satisfies the relationship
\[
c_s=\f{c_f}{1+\kappa_s/\kappa_c},
\]
where $\kappa_s$ is the surface reaction rate constant.  $\kappa$  and $\phi$
are the  permeability and porosity of the rock, respectively. And
the permeability and the porosity satisfy the Carman-Kozeny correlation 
\[
\f{\kappa}{\kappa_0}=\f{\phi}{\phi_0}\left(\f{\phi(1-\phi_0)}{\phi_0(1-\phi)}\right)^2,
\]
where $\phi_0$ and $\kappa_0$ are the initial porosity and permeability of the rock respectively. 
$\alpha$ is the dissolving power of the acid and $\rho_s$ is the density of the solid phase. By the porosity
and permeability, $a_v$ is shown as
\[
\f{a_v}{a_0}=\f{\phi}{\phi_0}\sqrt{\f{\kappa_0\phi}{\kappa\phi_0}}=\f{1-\phi}{1-\phi_0},
\]
where $a_0$ is the initial interfacial area.

As is well established, acid treatment of carbonate reservoirs is a widely employed technique for enhancing oil and gas well productivity. Wormholing refers to the phenomenon where elongated, channel-like structures form and propagate within subsurface formations as a result of acid injection into a supercritical acid dissolution system. Given its critical role in improving the productivity of oil and gas reservoirs, wormhole propagation has garnered significant research attention over recent decades. Theoretical investigations into numerical methods for modeling these phenomena possess broad practical applicability and hold considerable significance. In \cite{rui0},  combining with the method of characteristics, Rui and Li investigated the block-centered finite difference method for simulating incompressible wormhole propagation. They subsequently extended this approach to compressible wormhole propagation in \cite{rui1}. The first author of this article, along with coauthors, developed a mass-preserving characteristic mixed finite element method for incompressible wormhole propagation in \cite{zhang1}. Additionally, they employed this mass-conservative technique to formulate a characteristic splitting mixed finite element method for compressible wormhole propagation in \cite{zhang2}. Due to high velocity and non-uniform porosity, Kou et al. established a fully conservative mixed finite element method within the Darcy-Forchheimer framework for incompressible wormhole problems in \cite{sun1}. Furthermore, they explored a parallel algorithm for wormhole propagation under the Darcy-Brinkman-Forchheimer framework in \cite{sun2}. Generally, traditional mixed finite element methods or finite element methods for wormhole problems, as presented in \cite{zhang1,zhang2,sun1,sun2}, result in large-scale coupled systems, leading to higher computational costs. To address this issue, the first author of this article and coauthors introduced a hybrid mixed finite element method to simulate pressure and velocity in \cite{zhang2023,zhang2024}. By incorporating Lagrange multipliers at the edges of each element, the resulting global mixed system involves only the degrees of freedom associated with these multipliers, thereby significantly enhancing computational efficiency. Moreover, Guo, Tian, et al. examined the local discontinuous Galerkin finite element method for incompressible wormhole propagation in \cite{guo1}.

The primary objective of this article is to introduce a novel combined hybridized mixed discontinuous Galerkin finite element method for simulating incompressible wormhole propagation as described by Equation (1). Specifically, we employ a hybrid mixed finite element method for the pressure and velocity equations. Additionally, we extend this hybridization technique to formulate a hybrid mixed discontinuous Galerkin method for the convection-dominated concentration equation with an upwind scheme. To ensure the boundedness of porosity, which is crucial for both practical applications and theoretical analysis, we incorporate the ``cut-off" operator as introduced in Sun  \cite{sun1}. Compared to traditional discontinuous Galerkin methods, our proposed approach results in a global system with fewer unknowns and sparser stencils, thereby significantly reducing computational costs while effectively addressing discontinuous or fractured media problems. We rigorously establish the existence and uniqueness of the solution for the new combined method and derive optimal error estimates. Finally, numerical examples are provided to validate our theoretical findings and demonstrate the efficiency of the proposed method.

For convenience of analysis, we assume that $C$ and $\varepsilon$ are  some usual constant and small positive constant throughout this article,  which are independent of mesh parameter and time increment. 

\section{The formulation of combined method}

Denote $(\cdot,\cdot)$ the inner product in $L^{2}(\Omega)$ or $[L^{2}(\Omega)]^{d}$.  Introduce the divergence space: $H(\up{div};\Omega)=\{\mathbf{v}\in[L^{2}(\Omega)]^{d};~\nabla\cdot \mathbf{v}\in L^{2}(\Omega)\}$.  Next, we will formulate our method.

Set
\[
K=\f{a_{0}}{1-\phi_{0}}\f{\kappa_c\kappa_s}{\kappa_c+\kappa_s},\quad\beta(\phi)=\f{\mu}{\kappa(\phi)}.
\]
Then \eqref{e1} can be rewritten as follows
\begin{equation} \label{e2}
\begin{split}
&(\up{a})\quad\nabla\cdot\mathbf{u}+ \f{\alpha K(1-\phi)}{\rho_s}c_f=f, \\
&(\up{b})\quad \nabla p+\beta(\phi)\mathbf{u}=0,\\
&(\up{c})\quad \f{\p }{\p t}(\phi c_f)+\nabla\cdot(\mathbf{u}c_f-\phi D\nabla c_f)+K(1-\phi)c_f=f_Pc_f+f_Ic_I,\\
&(\up{d})\quad\f{\p \phi}{\p t}= \f{\alpha K(1-\phi)}{\rho_s}c_f.
\end{split}
\end{equation}

For the convenience  of analysis,  we  make the following assumptions:
\begin{assumption}\label{ass1}
Assume that the parameters
$\mu$, $k_c$, $k_s$, $\alpha$, $\rho_s$ are positive constants, and that $\phi_0$, $\frac{k(\phi)}{\mu}$ and $f(\cdot,t)$ are bounded as follows:
\begin{equation}\label{e3}
0<a_{\ast}\leq\frac{k(\phi)}{\mu}\leq a^{\ast},\quad 0<\phi_{0}< 1,\quad
| f(\cdot,t)|\leq C, 
\end{equation}
where  $a_{\ast}$, $a^{\ast}$ and $C$ are some positive constants. And we also assume that the diffusion coefficient $\mathbf{D}(\mathbf{u})$ satisfies the uniformly positive definiteness and Lipschitz continuousness 
\begin{equation}
\begin{aligned}\nonumber
&\mathbf{D}(\mathbf{u})\nabla c\cdot\nabla c\geq D_{*}| \nabla c|^2
\end{aligned}
\end{equation}
and 
\begin{equation}\label{e4}
\begin{aligned}
&\| \mathbf{D}(\mathbf{u})-\mathbf{D}(\mathbf{v})\|_{[L^2]^d}\leq D^{*}\| \mathbf{u}-\mathbf{v}\|_{[L^2]^d},
\end{aligned}
\end{equation}
where $D*$ and $D_{*}$ are two positive constants independent of $\mathbf{u}$ and $\mathbf{v}$ and $c$.
\end{assumption}

It is easily to obtain the following lemma.
\begin{lem}
Under the assumpution \ref{ass1}, we know that $\kappa^{-1}(\phi )$ is bounded and Lipschitz continuous, that is, there exists some constant $C$ such that
\[
\kappa^{-1}(\phi )\leq C,\quad |\kappa^{-1}(\phi_1)-\kappa^{-1}(\phi_1)|\leq C|\phi_1-\phi_2|.
\]
\end{lem}


To establish our method for system \eqref{e1},  let $\mathcal{T}_h$ be a shape regular partition of $\Omega$ with $\mathcal{T}_h = \{K_1, K_2, . . . , K_N \}$, and denote $\partial\mathcal{T}_h=\cup_{K\in\mathcal{T}_h}\{e|\,e\in\partial K\}$ to be the set of all cell edges. The velocity vector field $\mathbf{u}$ induces a natural splitting of element boundaries into inflow and outflow parts, i.e., we denote the outflow boundary $\partial K^{out}=\{x\in\partial K:\mathbf{u}\cdot\mathbf{n}>0\}$ and $\partial K^{in}=\partial K\backslash\partial K^{out}$, where $\mathbf{n}$ denotes the unit outward normal direction of $\p K$. The unions of the element inflow and outflow boundaries are denoted by $\partial\mathcal{T}_h^{in}$ and $\partial\mathcal{T}_h^{out}$, respectively. And $\partial\Omega^{in}$ and $\partial\Omega^{out}$ are the inflow and outflow regions of the boundary $\partial\Omega$. Furthermore, let $h_e = \textrm{diam}(e)$ and $h=\max_e(h_e)$ for all $e\in\partial\mathcal{T}_h$.

We need to introduce the piecewise Sobolev spaces. For any integer $s\geq 0$, we define
\begin{equation}\nonumber
H^s(\mathcal{T}_h)=\{v\in L^2(\Omega):v|_K\in H^s(K),K\in\mathcal{T}_h\},
\end{equation}
and
\begin{equation}\nonumber
\begin{aligned}
&L^2(\partial\mathcal{T}_h)=\{v\in L^2(e),\forall e\in\partial\mathcal{T}_h\}.
\end{aligned}
\end{equation}
Define some inner products 
\begin{equation}\nonumber
\begin{aligned}
&(u,v)_{K}=\int_{K}uvdx,\quad (u,v)_{\mathcal{T}_h}=\sum\limits_{K\in\mathcal{T}_h}(u,v)_K,\\
 & \langle u,v\rangle_{\p  K}=\int_{\p K}uvds,\quad \langle u,v\rangle_{\partial\mathcal{T}_h}=\sum\limits_{K\in\mathcal{T}_h}\langle u,v\rangle_{\p K},
\end{aligned}
\end{equation}
and the norms $\|\cdot\|_{\mathcal{T}_{h}}=\sqrt{(\cdot,\cdot)_{\mathcal{T}_h}}$ and  $|\cdot|_{\partial\mathcal{T}_{h}}=\sqrt{\langle\cdot,\cdot\rangle_{\partial\mathcal{T}_h}}$.

Introduce the discrete approximate spaces  $\Psi_h$, $\Lambda_h$, $\Pi_h$ and $\Sigma_h$
\begin{equation}\nonumber
\begin{aligned}
&\Psi_h=\{v\in H^k(\mathcal{T}_h):v|_K\in P_k(K),K\in\mathcal{T}_h\},\\
&\Lambda_h=\{v\in L^2(\mathcal{T}_h):v|_K\in P_k(K),K\in\mathcal{T}_h\},\\
&\Theta_h=\{v\in [L^2(\mathcal{T}_h)]^d:v|_K\in RT_k(K),K\in\mathcal{T}_h\},\\
&\Sigma_h=\{v\in L^2(\partial\mathcal{T}_h):v|_e\in P_k(e),e\in\partial\mathcal{T}_h\},
\end{aligned}
\end{equation}
where $P_k(K)$ and $P_k(e)$ are the spaces of polynomial functions of degree at most $k$ for each $K\in\mathcal{T}_h$ and each $e\in\partial\mathcal{T}_h$, respectively, and $RT_k(K)=[P_k(K)]^d\oplus xP_k(K)$ denotes the Raviart-Thomas element space as in \cite{BF1991}.

In order to give the time discrete formulation of system \eqref{e2}, we  let $N>0$ be a positive integer and $\Delta t=T/N$ be the time size, and  denote $t^{n}=n\Delta t$, $n=0,1,\cdots, N$. For any function $w$, denote $w^n$  the value of $w$ at time $t^n$.

Next, we will formulate our method for  incompressible wormhole propagation problem.

\subsection{The discretization of the porosity eqaution}

Here we consider the similar `` cut-off " technique for the discretization of the porosity as in \cite{zhang2023,sun1}, which is very important role both in practical production and our theoretical analysis.  We revisit numerical procedure of the porosity:
\begin{alg}\label{alg1}
Give an initial approximation
$(\phi_h^0,c_h^0)\in 
{\Lambda}_{h}\times\Psi_{h}$, 
for $n=1,2,\ldots, N$,   seek $\phi_h^n$ such that
\begin{equation}\label{e6}
\frac{\phi^n_h-\phi^{n-1}_h}{\Delta t}=\f{\alpha K (1-\phi^{n}_h)}{\rho_s}\bar{c}^{n-1}_h,\quad \bar{c}^{n-1}_h=\max(0,\min(c^{n-1}_{h},1))
\end{equation}
where  $c_{h}$ is a given approximation of the concentration $c_{f}$.
\end{alg}

\begin{thm}\label{thm1}
The discrete scheme \eqref{e6} preserves the maximum principle for any time step size, that is, if we suppose that $0<\phi^{0}_{h}<1$, we have that $0<\phi^{0}_{h}\leq\phi^{n}_{h}\leq\phi^{n+1}_{h}<1$ for any $n\geq 1$.
\end{thm}
\begin{proof}
From \eqref{e6}, we have
\[
\phi^n_h=\f{\phi^{n-1}_h}{1+\theta_{n}}+\f{\theta_{n}}{1+\theta_{n}},\quad \theta_{n}=\f{\Delta t\alpha K}{\rho_{s}}\bar{c}^{n-1}_{h}.
\]
Under the assumption $0<\phi^{0}_{h}<1$, using the recurrence relation, we can easily get the estimate
\begin{equation}\label{e7}
0<\phi^n_h<1,\quad n= 1,2,\cdots, N.
\end{equation}
By \eqref{e7} we get
\[
\phi^{n+1}_h-\phi^{n}_{h}=\f{\phi^{n}_h}{1+\theta_{n+1}}+\f{\theta_{n+1}}{1+\theta_{n+1}}-\phi^{n}_{h}=(1-\phi^{n}_{h})\f{\theta_{n+1}}{1+\theta_{n+1}}.
\]
By the recurrence relation and \eqref{e7}, we get the estimate
\[
0<\phi^{0}_{h}\leq \phi^{n}_h\leq \phi^{n+1}_{h}<1.
\]
\end{proof}

\subsection{Hybrid mixed finite element method for the pressure and Darcy velocity equations}
In this subsection, we give the hybrid mixed finite element  method for pressure and velocity as  in \cite{zhang2023}. Here we use completely discontinuous piecewise polynomial functions and ensure the continuity of the normal fluxes over internal interfaces by adding Lagrangian multiplier. The hybrid mixed finite element formulation can be written as below:
\begin{alg}\label{alg2}
For given approximate value  $c_{h}^{n-1}$, seek $(\mathbf{u}_h^{n}, p_h^{n}, \lambda_h^{n}) \in  \Theta_h \times \Lambda_h \times \Sigma_h$ such that
\begin{equation}\label{e8}
\begin{aligned}
(\nabla\cdot\mathbf{u}_h^{n},v_h)_{\mathcal{T}_h}+(\f{\alpha K(1-\phi^n_h)}{\rho_s}\bar{c}^{n-1}_h,v_h)_{\mathcal{T}_h}&=(f^{n},v_h)_{\mathcal{T}_h},\quad&\forall v_h \in \Lambda_h,\\
(\beta(\phi_{h}^{n})\mathbf{u}_h^{n},\bm\tau_h)_{\mathcal{T}_h}-(p_h^{n},\nabla\cdot\bm\tau_h)_{\mathcal{T}_h}+\langle\lambda_h^{n},\bm\tau_h\cdot \mathbf{n}\rangle_{\partial\mathcal{T}_h}&=0,\quad&\forall \bm\tau_h \in \Theta_h,\\
\langle\mathbf{u}_h^{n}\cdot \mathbf{n},\mu_h\rangle_{\partial\mathcal{T}_h}&=0,\quad&\forall \mu_h \in \Sigma_h,
\end{aligned}
\end{equation}
where $\beta(\phi_{h}^{n})=\kappa(\phi_{h}^{n})/\mu$.
\end{alg}

Define the symmetric bilinear form:
\[
\begin{aligned}
&B_{\mathbf{u}}(\phi^n_h;(\mathbf{u}_h^{n},p_h^{n},\lambda_h^{n}),(\bm\tau_h,v_h,\mu_h))\\&\quad:=(\beta(\phi_h^{n})\mathbf{u}_h^{n},\bm\tau_h)_{\mathcal{T}_h}+(\nabla p_h^{n},\bm\tau_h)_{\mathcal{T}_h}+(\mathbf{u}_h^{n},\nabla v_h)_{\mathcal{T}_h}\\
&\quad\quad+\langle\lambda_h^{n}-p_h^{n},\bm\tau_h\cdot \mathbf{n}\rangle_{\partial\mathcal{T}_h}+\langle\mathbf{u}_h^{n}\cdot \mathbf{n},\mu_h-v_h\rangle_{\partial\mathcal{T}_h}.
\end{aligned}
\]
We can rewrite \eqref{e8} into the following equivalent form:
\begin{alg}[HMFE Algorithm]\label{alg3}
 For given approximate values $\phi^n_h$ and  $c_h^{n-1}$, find $(\mathbf{u}_h^{n}, p_h^{n}, \lambda_h^{n,p}) \in \Theta_h \times \Lambda_h \times \Sigma_h$ such that
\[
B_{\mathbf{u}}(\phi^n_h;(\mathbf{u}_h^{n},p_h^{n},\lambda_h^{n}),(\bm\tau_h,v_h,\mu_h))=(\f{\alpha K(1-\phi^n_h)}{\rho_s}\bar{c}^{n-1}_h-f^{n},v_h)_{\mathcal{T}_h}, 
\]
\[\forall (\bm\tau_h, v_h, \mu_h) \in \Theta_h \times \Lambda_h \times \Sigma_h.
\]
\end{alg}

\subsection{Mixed hybridized  discontinuous Galerkin method for the concentration equation}
We can rewrite the concentration equation of \eqref{e2} into the  first order partial differential equations:
\begin{equation}\nonumber
\begin{aligned}
\bm\sigma+\phi D(\mathbf{u})\nabla c_f&=0,\\
\frac{\partial }{\partial t}(\phi c_f)+\nabla\cdot(\bm\sigma+\mathbf{u}c_f)+K(1-\phi)c_f&=f_Ic_I+f_Pc_f.
\end{aligned}
\end{equation}
Multiplying the above two equations by the test functions $\bm\tau$ and $v$ respectively,  and adding an upwind stabilization term, we can reach the mixed variational formulation for the concentration equation
\begin{equation}\label{e9}
\begin{aligned}
&([\phi D(\mathbf{u})]^{-1}\bm\sigma,\bm\tau)_{\mathcal{T}_h}+(\nabla c_f,\bm\tau)_{\mathcal{T}_h}=0,\\
&(\frac{\partial  }{\partial t}(\phi c_f),v)_{\mathcal{T}_h}+(\nabla\cdot(\bm\sigma+\mathbf{u}c_f),v)_{\mathcal{T}_h}+\langle\mathbf{u}\cdot\mathbf{n}(c^{+}-c_f),v\rangle_{\partial\mathcal{T}^{in}_h}
\\&\qquad\qquad\qquad\qquad\qquad+(K(1-\phi)c_f,v)_{\mathcal{T}_h}=(f_Ic_I+f_Pc_f,v)_{\mathcal{T}_h},
\end{aligned}
\end{equation}
where $c^{+}:=c_f|_{\partial K^{+}}$ denotes the upwind value at the upwind element $K^+$, that is, the element attached to $e$ where $\mathbf{u}\cdot\mathbf{n}=\mathbf{u}\cdot\mathbf{n}_{e}\geq 0$. To incorporate the boundary condition, we define $c^+=0$
on $\partial\Omega^{in}$. After integration by parts, we can reach
\begin{equation}\nonumber
\begin{aligned}
&([\phi D(\mathbf{u})]^{-1}\bm\sigma,\bm\tau)_{\mathcal{T}_h}+(\nabla c_f,\bm\tau)_{\mathcal{T}_h}=0,\\
& (\frac{\partial  }{\partial t}(\phi c_f),v)_{\mathcal{T}_h}-(\bm\sigma+\mathbf{u}c_f,\nabla v)_{\mathcal{T}_h}+\langle\mathbf{u}\cdot\mathbf{n}c^{+},v\rangle_{\partial\mathcal{T}_h^{in}}+\langle\mathbf{u}\cdot\mathbf{n}c_f,v\rangle_{\partial\mathcal{T}_h^{out}}
\\&\qquad+(K(1-\phi)c_f,v)_{\mathcal{T}_h}=(f_Ic_I+f_Pc_f,v)_{\mathcal{T}_h}.
\end{aligned}
\end{equation}

Introduce the upwind value as a new variable $\lambda^{c}:=c^+$, and define $\hat{c}$ as follows: for any $K\in\mathcal{T}_{h}$
\begin{equation}\label{e10}
\hat{c}:=\left\{
\begin{aligned}
&\lambda^c,& \quad e\subset\partial K^{in},\\
&c_f,& \quad e\subset\partial K^{out}.\\
\end{aligned}
\right.
\end{equation}
Using the fact that $\lambda^{c}=\hat{c}=c^+$ on both sides of $e$, we can give the corresponding fully discrete  mixed hybridized discontinuous Galerkin method for the concentration with the first-oder backward Euler difference scheme in time.
 \begin{alg}\label{alg4}
 For given approximate values $\phi^{n}_h$, $c^{n-1}_h$ and $\mathbf{u}_{h}^{n}$, seek $(\bm\sigma_h^{n}, c_h^{n}, \lambda_h^{c,n}) \in  \Theta_h \times \Lambda_h \times \Sigma_h$
\begin{equation}\label{e11}
\begin{aligned}
&([\phi^n_hD(\mathbf{u}^{n}_h)]^{-1}\bm\sigma^{n}_h,\bm\tau_h)_{\mathcal{T}_h}-(c_h^{n},\nabla\cdot\bm\tau_h)_{\mathcal{T}_h}+\langle\lambda^{c,n}_h,\bm\tau_h\cdot\mathbf{n}\rangle_{\partial\mathcal{T}_h}=0,\quad\bm\tau_h\in\Theta_h,\\
\end{aligned}
\end{equation}
\begin{equation}\label{e12}
\begin{aligned}
(\f{\phi^n_hc_h^{n}-\phi_h^{n-1}c_h^{n-1}}{\Delta t},v_h)_{\mathcal{T}_h}&+(\nabla\cdot(\mathbf{u}^{n}_hc_h^{n})+\nabla\cdot\bm\sigma_h^{n},v_h)_{\mathcal{T}_h}+\langle\mathbf{u}_h^{n}\cdot\mathbf{n}(\lambda_h^{c,n}-c_h^{n}),v_h\rangle_{\partial\mathcal{T}_h^{in}}
\\
&+(K(1-\phi^n_h)c^n_h-f^n_Pc^n_h,v_h)_{\mathcal{T}_h}=(f^n_Ic^n_I,v_h)_{\mathcal{T}_h},\quad v_h\in\Lambda_h,\\
\end{aligned}
\end{equation}
\begin{equation}\label{e13}
\langle\bm\sigma^{n}_h\cdot\mathbf{n},\mu_h\rangle_{\partial\mathcal{T}_h}=0,\quad\mu_h\in\Sigma_h,
\end{equation}
\begin{equation}\label{e14}
\langle\mathbf{u}^{n}_h\cdot\mathbf{n}\hat{c}_h^{n},\mu_h\rangle_{\partial\mathcal{T}_h}=0,\quad\mu_h\in\Sigma_h.
\end{equation}
\end{alg}

Define the bilinear form:
\[
\begin{aligned}
&B_c((\phi^n_h,\mathbf{u}^n_h);(\bm\sigma_h^{n},c_h^{n},\lambda_h^{c,n}),(\bm\tau_h,v_h,\mu_h))
\\&\quad:=-\frac{1}{\Delta t}(\phi^n_h c_h^{n},v_h)_{\mathcal{T}_h}-([K(1-\phi^n_h)-f^n_P]c^n_h,v_h)_{\mathcal{T}_h}+(\bm\sigma_h^{n}+\mathbf{u}^{n}_{h}c_h^{n},\nabla v_h)_{\mathcal{T}_h}\\&\quad\quad+([\phi^n_hD(\mathbf{u}_h^{n})]^{-1}\bm\sigma_h^{n}+\nabla c_h^{n},\bm\tau_h)_{\mathcal{T}_h}
+\langle\lambda_h^{c,n}-c_h^{n},\bm\tau_h\cdot\mathbf{n}\rangle_{\partial\mathcal{T}_h}\\
&\quad\quad+\langle\bm\sigma_h^{n}\cdot\mathbf{n}+\mathbf{u}_h^{n}\cdot\mathbf{n}\hat{c}_h^{n},\mu_h-v_h\rangle_{\partial\mathcal{T}_h}.
\end{aligned}
\]
Now, we can arrive at the following upwind  mixed hybridized discontinuous Galerkin finite element (UMHDG) method.
\begin{alg}[UMHDG Algorithm]\label{alg5}
 For given approximate values $\phi^n_h$, $c^{n-1}_h$ and $\mathbf{u}_h^{n}$, find $(\bm\sigma_h^{n},c_h^{n},\lambda_h^{c,n})\in \Theta_h \times \Lambda_h \times \Sigma_h$ such that
\[
B_c((\phi^n_h,\mathbf{u}^n_h);(\bm\sigma_h^{n},c_h^{n},\lambda_h^{c,n}),(\bm\tau_h,v_h,\mu_h))=-(f^{n}_Ic^{n}_{I}+\phi^{n-1}_h c_h^{n-1}/ \Delta t,v_h)_{\mathcal{T}_h},\]
\[ \forall (\bm\tau_h,v_h,\mu_h) \in \Theta_h \times \Lambda_h \times \Sigma_h.
\]
\end{alg}

\subsection{The combined  mixed hybridized discontinuous Galerkin finite element method}

Now we give the combined hybrid mixed discontinuous Galerkin finite element method for incompressible wormhole propagation with upwind technique.
\begin{alg}[Combined MHDG Algorithm]\label{alg6}
For the given initial approximate values $(\phi_h^0,c_h^0)\in 
{\Lambda}_{h}\times\Psi_{h}$, 
for $n=1,2,\ldots, N$,   seek $\phi_h^n$, $(\mathbf{u}_h^{n}, p_h^{n}, \lambda_h^{n}) \in \Theta_h \times \Lambda_h \times \Sigma_h$ and $(\bm\sigma_h^{n},c_h^{n},\lambda_h^{c,n})\in \Theta_h \times \Lambda_h \times \Sigma_h$, such that, for $ \forall (\bm\tau_h,v_h,\mu_h) \in \Theta_h \times \Lambda_h \times \Sigma_h$
\begin{equation}\label{e15}
\begin{aligned}
&({\rm a})\quad\frac{\phi^n_h-\phi^{n-1}_h}{\Delta t}=\f{\alpha K (1-\phi^{n}_h)}{\rho_s}\bar{c}^{n-1}_h,\\
&({\rm b})\quad B_{\mathbf{u}}(\phi^n_h;(\mathbf{u}_h^{n},p_h^{n},\lambda_h^{n}),(\bm\tau_h,v_h,\mu_h))=(\f{\alpha K(1-\phi^n_h)}{\rho_s}\bar{c}^{n-1}_h-f^{n},v_h)_{\mathcal{T}_h}, \\
&({\rm c})\quad B_c((\phi^n_h,\mathbf{u}^n_h);(\bm\sigma_h^{n},c_h^{n},\lambda_h^{c,n}),(\bm\tau_h,v_h,\mu_h))=-(f^{n}_Ic^{n}_{I}+\f{\phi^{n-1}_h c_h^{n-1}}{\Delta t},v_h)_{\mathcal{T}_h},
\end{aligned}
\end{equation}
where $\bar{c}^{n-1}_h=\max(0,\min(c^{n-1}_{h},1)).$
\end{alg}

\subsection{Existence and uniqueness}
In order to prove the existence and uniqueness of the solution of our proposed algorithm, we will use the following important result (see Lemma 3.1 in \cite{ES2010}).
\begin{lem}\label{lem1}
There is a unique solution $\tilde{\tau}_c\in \Theta_h$ defined element-wise by the variational problem
\[
\begin{aligned}
(\tilde{\tau}_c,\omega)_K&=(\nabla c_h,\omega)_K,\quad\quad\forall \omega\in[P_{k-1}(K)]^d,\\
\langle\tilde{\tau}_c\cdot\mathbf{n},\mu\rangle_{\partial K}&=\langle\frac{1}{h}(\lambda^{c}_h-c_h),\mu\rangle_{\partial K},\quad\quad \mu\in P_k(\partial K),
\end{aligned}
\]
where  $c_h\in\Lambda_h$ and $\lambda^{c}_h\in\Sigma_h$ are any piecewise polynomial functions.

Moreover, the following result
\[
\begin{aligned}
\|\tilde{\tau}_c\|_{\mathcal{T}_h}\leq C_c(\|\nabla c_h\|_{\mathcal{T}_h}^2+\frac{1}{h}|\lambda^{c}_h-c_h|_{\partial\mathcal{T}_h}^2)^{\frac{1}{2}}
\end{aligned}
\]
holds, where  $C_c$  is  a constant independent of the parameter $h$.
\end{lem}

Define the pair of the norms
\[
\begin{aligned}
&\|(\bm\omega,v,\mu)\|_\mathbf{u}:=(\|\bm\omega\|_{\mathcal{T}_h}^2+\|\nabla v\|_{\mathcal{T}_h}^2+\frac{1}{h}|\mu-v|_{\partial\mathcal{T}_h}^2)^\frac{1}{2},\\
&\|(\bm\omega,v,\mu)\|_\mathbf{u,*}:=(\|(\bm\omega,v,\mu)\|^{2}_\mathbf{u}+h|\bm\omega\cdot\mathbf{n}|^{2}_{\p \mathcal{T}_{h}})^\frac{1}{2}.
\end{aligned}
\]
Using the similar techniques as in \cite{zhang2023,zhang2024,zhang2017,zhang2020}, we can easily get the following stability and boundedness of the bilinear form $B_\mathbf{u}$.
\begin{lem}[Stability and boundedness of $B_{\mathbf{u}}$]\label{lem2}
There exist two positive constants $C^{\ast}$ and $C_{\ast}$,
which are independent of the mesh size $h$, such that
\begin{equation}\nonumber
\begin{aligned}
&(\up{a})\quad| B_{\mathbf{u}}(\phi;(\mathbf{u}, p, \lambda),(\bm\tau_h, v_h, \mu_h))|\leq C^{\ast}\|(\mathbf{u}, p, \lambda)\|_{\mathbf{u,*}}\|(\bm\tau_h, v_h, \mu_h)\|_\mathbf{u},\\
&(\up{b})\quad\sup\limits_{(\bm\tau_h, v_h, \mu_h)}\frac{B_{\mathbf{u}}(c_h;(\mathbf{u}_h, p_h, \lambda_h),(\bm\tau_h, v_h, \mu_h))}{\|(\bm\tau_h, v_h, \mu_h)\|_\mathbf{u}}\geq C_{\ast}\|(\mathbf{u}_h, p_h, \lambda_h)\|_\mathbf{u},
\end{aligned}
\end{equation}
for all $( \mathbf{u}, p, \lambda) \in  H^k(\mathcal{T}_h) \times L^2(\mathcal{T}_h)\times L^2
(\partial\mathcal{T}_h)$, $( \mathbf{u}_h, p_h, \lambda_h) \in \Theta_h \times \Lambda_h \times \Sigma_h$ and $(\bm\tau_h, v_h, \mu_h) \in \Theta_h \times \Lambda_h \times \Sigma_h$.
\end{lem}

To show the stability and boundedness of the bilinear form $B_c$, we define the norms as follows
\[
\begin{aligned}
\|(\bm\sigma,c,\lambda^c)\|_{D}&:=(\|\bm\sigma\|_{\mathcal{T}_h}^2+\f{1}{\Delta t}\| c\|_{\mathcal{T}_h}^2+\|\nabla c\|_{\mathcal{T}_h}^2+\frac{1}{h}|\lambda^c-c|_{\partial\mathcal{T}_h}^2)^\frac{1}{2},\\
\|(c,\lambda^c)\|_B&:=(\|\mathbf{u}\cdot\nabla c\|_{\mathcal{T}_h}^2+\|\mathbf{u}\cdot\mathbf{n}|^{1/2}(\lambda^c-c)|_{\partial\mathcal{T}_h}^2)^\frac{1}{2},\\
\|(\bm\sigma,c,\lambda^c)\|_{B,*}&:=(h|\bm\sigma\cdot\mathbf{n}|_{\p\mathcal{T}_h}^2+\|\mathbf{u}\cdot\nabla c\|_{\partial\mathcal{T}_h}^2+|\mathbf{u}\cdot\mathbf{n}\lambda^{c}|^{2}_{\p\mathcal{T}_{h}})^\frac{1}{2},\\
\|(\bm\sigma,c,\lambda^c)\|_c&:=(\|(\bm\sigma,c,\lambda^c)\|_D^2+\|(c,\lambda^c)\|_B^2)^{\frac{1}{2}},
\\
\|(\bm\sigma,c,\lambda^c)\|_{c,*}&:=(\|(\bm\sigma,c,\lambda^c)\|_D^2+\|(c,\lambda^c)\|_{B,*}^2)^{\frac{1}{2}}.
\end{aligned}
\]
We have the following the stability and boundedness result on  the bilinear form $B_c$.

\begin{lem}[Stability and boundedness of $B_{c}$]\label{lem3}
There exist  two positive constants $C^{c,\ast}$,  and $C_{c,\ast}$, which are
independent of the mesh size $h$, such that, for some given $\Delta t_{0}>0$, when $\Delta t\leq \Delta t_{0}$,
\begin{equation}\nonumber
\begin{aligned}
&(\up{a})\quad| B_c((\phi,\mathbf{u});(\bm\sigma,c,\lambda^{c}),(\bm\tau_h,v_h,\mu_h))|\leq C^{c,\ast}\|(\bm\sigma,c,\lambda^{c})\|_{c,*}\|(\bm\tau_h,v_h,\mu_h)\|_c,\\
&(\up{b})\quad\sup\limits_{(\bm\tau_h,v_h,\mu_h)}\frac{B_c((\phi_h,\mathbf{u}_h);(\bm\sigma_h,c_h,\lambda_h^{c}),(\bm\tau_h,v_h,\mu_h))}{\|(\bm\tau_h, v_h, \mu_h)\|_c}\geq C_{c,\ast}\|(\bm\sigma_h,c_h,\lambda_h^{c})\|_c,\\
&(\up{c})\quad\forall ( \bm\sigma_h,c_h,\lambda_h^{c}) \in \Theta_h \times \Lambda_h \times \Sigma_h,\quad B_c((\phi_h,\mathbf{u}_h);(\bm\sigma_h,c_h,\lambda_h^{c}),(\bm\tau_h,v_h,\mu_h))=0  \Rightarrow(\bm\tau_h, v_h, \mu_h)= 0,
\end{aligned}
\end{equation}
holds for all $( \bm\sigma_h,c_h,\lambda_h^{c}) \in \Theta_h \times \Lambda_h \times \Sigma_h$ and $(\bm\tau_h, v_h, \mu_h) \in \Theta_h \times \Lambda_h \times \Sigma_h$.
\end{lem}

\begin{proof}
Firstly, we choose $(\bm\tau_{h},v_{h},\mu_{h})=(\gamma\tilde{\bm\tau}_{c},0,0)$ in the bilinear form $B_{c}$ and use Lemma \ref{lem1} to get
\begin{equation}\label{e16}
\begin{split}
&B_{c}((\phi_h,\mathbf{u}_h);(\bm\sigma_h,c_h,\lambda_h^{c}),(\gamma\tilde{\bm\tau}_{c},0,0))
\\=&(\frac{1}{\phi_hD(\mathbf{u}_h)}\bm\sigma_h,\gamma\tilde{\bm\tau}_{c})_{\mathcal{T}_h}+(\nabla c_{h},\gamma\tilde{\bm\tau}_{c})_{\mathcal{T}_h}+\langle\lambda_h^{c}-c_h,\gamma\tilde{\bm\tau}_{c}\cdot\mathbf{n}\rangle_{\partial\mathcal{T}_h}
\\=&\gamma(\frac{1}{\phi_hD(\mathbf{u}_h)}\bm\sigma_h,\tilde{\bm\tau}_{c})_{\mathcal{T}_h}+\gamma\|\nabla c_h\|^{2}_{\mathcal{T}_h}+\f{\gamma}{h}|\lambda_h^{c}-c_h|^{2}_{\partial\mathcal{T}_h}
\\\geq&-\f{1}{2}(\frac{1}{\phi_hD(\mathbf{u}_h)}\bm\sigma_h,\bm\sigma_h)_{\mathcal{T}_h}-\f{\gamma^{2}}{2\phi^0_hD_*}\|\tilde{\bm\tau}_{c}\|^{2}_{\mathcal{T}_h}+\gamma(\|\nabla c_h\|^{2}_{\mathcal{T}_h}+\f{1}{h}|\lambda_h^{c}-c_h|^{2}_{\partial\mathcal{T}_h})
\\\geq&-\f{1}{2}(\frac{1}{\phi_hD(\mathbf{u}_h)}\bm\sigma_h,\bm\sigma_h)_{\mathcal{T}_h}+(\gamma-\f{C^{2}_{c}\gamma^{2}}{2\phi^0_hD_*})(\|\nabla c_h\|^{2}_{\mathcal{T}_h}+\f{1}{h}|\lambda_h^{c}-c_h|^{2}_{\partial\mathcal{T}_h}).
\end{split}
\end{equation}
And then, taking $(\bm\tau_{h},v_{h},\mu_{h})=(\bm\sigma_{h},-c_{h},-\lambda^{c}_{h})$  in the bilinear form $B_{c}$, we have
\begin{equation}\label{e17}
\begin{split}
&B_{c}((\phi_h,\mathbf{u}_{h});(\bm\sigma_h,c_h,\lambda_h^{c}),(\bm\sigma_{h},-c_{h},-\lambda^{c}_{h}))
\\=&\f{1}{\Delta t}(\phi_h c_h,c_h)_{\mathcal{T}_h}+([K(1-\phi_h)-f_P]c_h,c_h)_{\mathcal{T}_h}+(\frac{1}{\phi_hD(\mathbf{u}_h)}\bm\sigma_h,\bm\sigma_h)_{\mathcal{T}_h}
\\&-(\mathbf{u}_hc_h,\nabla c_h)_{\mathcal{T}_h}+\langle\mathbf{u}_h\cdot\mathbf{n}\hat{c}_h,-\lambda_h^{c}+c_h\rangle_{\partial\mathcal{T}_h}
\\=&\f{1}{\Delta t}(\phi c_h,c_h)_{\mathcal{T}_h}+([K(1-\phi_h)-f_P]c_h,c_h)_{\mathcal{T}_h}+\f{1}{2}(\nabla\cdot\mathbf{u}_{h} c_h,c_h)_{\mathcal{T}_h}
\\&+(\frac{1}{\phi_hD(\mathbf{u}_h)}\bm\sigma_h,\bm\sigma_h)_{\mathcal{T}_h}
-\f{1}{2}\langle\mathbf{u}_h\cdot\mathbf{n}{c}_h,c_h\rangle_{\partial\mathcal{T}_h}+\langle\mathbf{u}_h\cdot\mathbf{n}\hat{c}_h,-\lambda_h^{c}+c_h\rangle_{\partial\mathcal{T}_h}
\\=&\f{1}{\Delta t}(\phi c_h,c_h)_{\mathcal{T}_h}+([K(1-\phi_h)-f_P]c_h,c_h)_{\mathcal{T}_h}+(\frac{1}{\phi_hD(\mathbf{u}_h)}\bm\sigma_h,\bm\sigma_h)_{\mathcal{T}_h}
\\&-\f{1}{2}\langle\mathbf{u}_h\cdot\mathbf{n}{c}_h,c_h\rangle_{\partial\mathcal{T}_h}+\langle\mathbf{u}_h\cdot\mathbf{n}\hat{c}_h,-\lambda_h^{c}+c_h\rangle_{\partial\mathcal{T}_h}.
\end{split}
\end{equation}
Note that
\[
\begin{aligned}
\mathbf{u}_h\cdot\mathbf{n}\large |_{\partial\mathcal{T}^{in}_h}=-|\mathbf{u}_h\cdot\mathbf{n}| \large |_{\partial\mathcal{T}^{in}_h}, \mathbf{u}_h\cdot\mathbf{n}\large |_{\partial\mathcal{T}^{out}_h}=|\mathbf{u}_h\cdot\mathbf{n}|\large |_{\partial\mathcal{T}^{out}_h}, \hat{c}_h\large |_{\partial\mathcal{T}^{in}_h}=\lambda^c_h\large |_{\partial\mathcal{T}^{in}_h},
\hat{c}_h\large |_{\partial\mathcal{T}^{out}_h}=c_h\large |_{\partial\mathcal{T}^{out}_h}.
\end{aligned}
\]
So we have
\begin{equation}\label{e18}
\begin{aligned}
&(\up{a})\quad-\frac{1}{2}\langle\mathbf{u}_h\cdot\mathbf{n}c_h,c_h\rangle_{\partial\mathcal{T}_h}=\frac{1}{2}\langle|\mathbf{u}_h\cdot\mathbf{n}|c_h,c_{h}\rangle_{\partial\mathcal{T}_h^{in}}-\frac{1}{2}\langle|\mathbf{u}_h\cdot\mathbf{n}|c_h,c_{h}\rangle_{\partial\mathcal{T}_h^{out}},\\
&(\up{b})\quad-\langle\mathbf{u}_h\cdot\mathbf{n}\hat{c}_h,\lambda_h^{c}\rangle_{\partial\mathcal{T}_h}=\langle|\mathbf{u}_h\cdot\mathbf{n}|\lambda_h^{c},\lambda_h^{c}\rangle_{\partial\mathcal{T}_h^{in}}-\langle|\mathbf{u}_h\cdot\mathbf{n}| \lambda_h^{c},c_h\rangle_{\partial\mathcal{T}_h^{out}},\\
&(\up{c})\quad\langle\mathbf{u}_h\cdot\mathbf{n}\hat{c}_h,c_h\rangle_{\partial\mathcal{T}_h}=\langle|\mathbf{u}_h\cdot\mathbf{n}|c_h,c_h\rangle_{\partial\mathcal{T}_h^{out}}-\langle |\mathbf{u}_h\cdot\mathbf{n}|\lambda_h^{c},c_h\rangle_{\partial\mathcal{T}_h^{in}}.
\end{aligned}
\end{equation}
Now let $K_1$ and $K_2$ denote two elements sharing the facet  $e=\partial K_1^{out}\bigcap \p K_2^{in}$. Since $\lambda_h^{c}$ is a single value function on $e$, we have $\lambda_h^{c}|_{\partial K_1^{out}}=\lambda_h^{c}|_{\partial K_2^{in}}$, which means that we can shift the terms only involving the Lagrange multiplier between neighboring elements. Hence,  summing \eqref{e18} up, we can  rewrite the last two terms of \eqref{e17} as follows
\[
\frac{1}{2}\langle |\mathbf{u}_h\cdot\mathbf{n}|(\lambda_h^{c}-c_{h}),\lambda_h^{c}-c_h\rangle_{\partial\mathcal{T}_h}.
\]
So we get
\begin{equation}\label{e19}
\begin{aligned}
&B_{c}((\phi_h,\mathbf{u}_{h});(\bm\sigma_h,c_h,\lambda_h^{c}),(\bm\sigma_{h},-c_{h},-\lambda^{c}_{h}))
\\=&\f{1}{\Delta t}(\phi_h c_h,c_h)_{\mathcal{T}_h}+([K(1-\phi_h)-f_P]c_h,c_h)_{\mathcal{T}_h}
\\&+(\frac{1}{\phi_hD(\mathbf{u}_h)}\bm\sigma_h,\bm\sigma_h)_{\mathcal{T}_h}+\frac{1}{2}\langle |\mathbf{u}_h\cdot\mathbf{n}|(\lambda_h^{c}-c_{h}),\lambda_h^{c}-c_h\rangle_{\partial\mathcal{T}_h}.
\end{aligned}
\end{equation}
In \eqref{e19}, we choose some time step $\Delta t_{0}>0$ such that 
$\f{\phi^0_h}{\Delta t}-f_P\geq  C_{*}>0$ when $\Delta t\leq \Delta t_{0}$, where $C_*$ is a constant independent of $h$ and $\Delta t$.  And we take $\gamma =\phi^0_hD_{*}/C^{2}_{c}$ ($\gamma-\f{C^{2}_{c}\gamma^{2}}{2\phi^0_hD_*}=\f{\phi^0_hD_{*}}{2C^{2}_{c}}>0$) in \eqref{e16}. Then combining \eqref{e16} and \eqref{e19}, we obtain the stability of the bilinear form $B_{c}$.  

Using the similar technique as above, we choose $( \bm\sigma_h,c_h,\lambda_h^{c}) = (\bm\tau_{h},-v_{h},-\mu_{h})$ in the bilinear form $B_c$ to reach
\begin{equation}\label{e20}
\begin{split}
&B_{c}((\phi_h,\mathbf{u}_{h});(\bm\tau_{h},-v_{h},-\mu_{h}),(\bm\tau_{h},v_{h},\mu_{h}))
\\=&\f{1}{\Delta t}(\phi_h v_h,v_h)_{\mathcal{T}_h}+([K(1-\phi_h)-f_P]v_h,v_h)_{\mathcal{T}_h}+(\frac{1}{\phi_hD(\mathbf{u}_h)}\bm\tau_h,\bm\tau_h)_{\mathcal{T}_h}
\\&-\f{1}{2}\langle\mathbf{u}_h\cdot\mathbf{n}{v}_h,v_h\rangle_{\partial\mathcal{T}_h}+\langle\mathbf{u}_h\cdot\mathbf{n}\hat{v}_h,-\mu_h+v_h\rangle_{\partial\mathcal{T}_h}
\\=&(\f{\phi_h}{\Delta t} v_h,v_h)_{\mathcal{T}_h}+([K(1-\phi_h)-f_P]v_h,v_h)_{\mathcal{T}_h}+(\frac{1}{\phi_hD(\mathbf{u}_h)}\bm\tau_h,\bm\tau_h)_{\mathcal{T}_h}
\\&+\frac{1}{2}\langle |\mathbf{u}_h\cdot\mathbf{n}|(\mu_h-v_{h}),\mu_h-v_h\rangle_{\partial\mathcal{T}_h}.
\end{split}
\end{equation}
Hence, for $\Delta t\leq \Delta t_{0}$, if $B_{c}((\phi_h,\mathbf{u}_{h});(\bm\tau_{h},-v_{h},-\mu_{h}),(\bm\tau_{h},v_{h},\mu_{h}))=0$, then $(\bm\tau_{h},v_{h},\mu_{h})=0$, which implies that the third conclusion of Lemma \ref{lem3} holds. 

Using Cauchy inequality we can easily get the boundedness of the bilinear form $B_{c}$. 

\end{proof}

\begin{thm}[Existence and Uniqueness]\label{thm2}
For given initial approximate values $\phi^{0}_{h}$ and $c^{0}_{h}$, there exists a parameter $\Delta t_{0}>0$, such that, when $\Delta  t\leq \Delta t_{0}$,  Algorithm \ref{alg6} exists a unique solution.
\end{thm}

\begin{proof}
By the stability and boundedness of the bilinear forms $B_{\mathbf{u}}$ and $B_{c}$, 
with the Banach-Ne\v{c}as-Babu\v{s}ka (BNB) theorem  (see \cite{santo,Zhu1995,QRZ1998,GQZ2000,ZGL2003} or Theorem 2.6 in \cite{eg2004} ) or the Babu\v{s}ka-Lax-Milgram (BLM) theorem in \cite{IR1989}, we know that Algorithm \ref{alg6} has a unique solution.
\end{proof}

\section{Convergence analysis}
Next, we will give some important projection operators and approximate properties, which will be used to show the convergence theorem of our proposed method.

Introduce the local $L^2$-projection operators
$\Pi_h$ and $\Pi_e$ as follows:
\[
\begin{aligned}
&(p-\Pi_hp,v_h)_K=0,\quad \forall v_h\in P_k(K),\\
&\langle\lambda-\Pi_e\lambda,\mu_h\rangle_e=0,\quad \forall \mu_h\in P_k(e),
\end{aligned}
\]
 where $K\in\mathcal{T}_h$, $e\in\partial\mathcal{T}_h$, $p\in L^2(K)$ and $\lambda\in L^2(e)$. These projection operators satisfy the following error estimates (see Lemma 3.9 and Lemma 3.10 in \cite{ES2010}, or  Theorem 4.4.20 in \cite{BS2002}).
\begin{lem}\label{lem4}
Suppose that  $p$ satisfies the regularity assumption $p\in{L^\infty(0,T;H^{k+1}(K))}$. Then, for the local $L^2$-projection operators $\Pi_h$ and $\Pi_e$, we have the estimate
\[
\begin{aligned}
&\|p-\Pi_hp\|_K\leq Ch^s\|p\|_{s,K},\quad 0\leq s\leq k+1,\\
&\|\nabla(p-\Pi_hp)\|_K\leq Ch^s\|p\|_{s+1,K},\quad 0\leq s\leq k,\\
&|p-\Pi_hp|_e+|p-\Pi_ep|_e\leq Ch^{s+\frac{1}{2}}\|p\|_{s+1,K},\quad 0\leq s\leq k,
\end{aligned}
\]
 where $C$ is a constant independent of $h$.
\end{lem}
\begin{lem}\label{lem5}
Suppose that $c$ satisfies the regularity assumption $c\in{L^\infty(0,T;H^{k+1}(K))}$.  Then, for any element $K$, we can reach
\begin{equation}\nonumber
\begin{aligned}
&\|\nabla(c-\Pi_hc)\|_K\leq Ch^s\|c\|_{s+1,K},\quad 0\leq s\leq k+1,\\
&|c-\Pi_ec|_{\partial K}\leq Ch^{s}|c|_{s,\partial K},\quad 0\leq s\leq k+1,
\end{aligned}
\end{equation}
where $C$ is a constant independent of  $h$.
\end{lem}
Similarly, the projection operators for functions on $\mathcal{T}_h$ and $\partial\mathcal{T}_h$ are defined element-wise and are denoted by the same symbols.
For $\bm{\omega}\in H(\up{div},K)$ we utilize the Raviart-Thomas interpolation projection \cite{BF1991,ES2010} defined by
\[
\begin{aligned}
&(\bm{\omega}-\Pi^{RT}\bm{\omega},\bm\tau_h)_K=0,\quad\forall\bm\tau_h\in[P_{k-1}(K)]^d,\\
&\langle(\bm{\omega}-\Pi^{RT}\bm{\omega})\cdot \mathbf{n}_e,\mu_h\rangle_e=0,\quad\forall\mu_h\in P_{k}(e),\quad e\in\partial K.
\end{aligned}
\]
We can reach the error estimate as follows:
\begin{lem}\label{lem6}
For the projection operator $\Pi^{RT}$ defined as above, we have the estimate
\[
\begin{aligned}
&\|\bm{\omega}-\Pi^{RT}\bm{\omega}\|_K+h^{\frac{1}{2}}|\bm{\omega}-\Pi^{RT}\bm{\omega}|_{\partial K}\leq Ch^s\|\bm{\omega}\|_{s,K},\quad \frac{1}{2}\leq s\leq k+1,\\
&\|\nabla\cdot(\bm{\omega}-\Pi^{RT}\bm{\omega})\|_K\leq Ch^s\|\nabla\cdot\bm{\omega}\|_{s,K},\quad 1\leq s\leq k+1,
\end{aligned}
\]
where $C$ is a constant independent of $h$.
\end{lem}

The following trace inequalities (see Section 2.1 in \cite{br2008}) will be also used to prove the convergence theorem.
\begin{lem}\label{lem7}
For $\forall v\in H^1(K)$, the trace inequalities are shown below
\[
\begin{aligned}
&\| v\|_{0,e}^2\leq C(h_e^{-1}\| v\|_{0,K}^2+h_e\| v\|_{1,K}^2),\\
&\| \nabla v\cdot \mathbf{n}_e\|_{0,e}^2\leq C(h_e^{-1}\| \nabla v\|_{0,K}^2+h_e\|\nabla^2v\|_{0,K}^2).
\end{aligned}
\]
\end{lem}

\subsection{The error estimate of the porosity}

Set 
\[
\phi^n_h-\phi^n=e^n_\phi,\quad 
c^n_h-c^n_f=(c^n_h-\Pi_hc^n_f)+(\Pi_hc^n_f-c^n_f)=\xi^n_c-\eta^n_c.
\]
Now we firstly estimate the boundedness of $e^{n}_{\phi}$. From \eqref{e1}(a) and \eqref{e15}(a), we get the residual equation
\begin{equation}\label{e21}
\begin{split}
( \f{e_{\phi}^{n}-e_{\phi}^{n-1}}{\Delta t},q)_{\mathcal{T}_h}=&(\f{\alpha K}{\rho_{s}}(\bar{c}_{h}^{n-1}-{c}_{f}^{n-1}),q)_{\mathcal{T}_h}-(\f{\alpha K}{\rho_{s}}(\phi^{n}_{h}\bar{c}_{h}^{n-1}-\phi^{n}{c}^{n-1}_{f}),q)_{\mathcal{T}_h}\\
 &+(\f{\p \phi}{\p t}-\f{{\phi}^{n}-{\phi}^{n-1}}{\Delta t}+\f{\alpha K(1-\phi^{n})}{\rho_{s}}({c}^{n-1}_{f}-c^{n}_{f}),q)_{\mathcal{T}_h},\quad q\in L^{2}(\Omega).
 \end{split}
\end{equation}

\begin{lem}\label{lem8}
There exists the following estimate
\begin{equation}\label{e22}
( \f{e_{\phi}^{n}-e_{\phi}^{n-1}}{\Delta t},q)_{\mathcal{T}_h}\leq C\{\|\xi_{c}^{n-1}\|^{2}_{\mathcal{T}_h}+\|\eta_{c}^{n-1}\|^{2}_{{\mathcal{T}_h}}+\|q\|^{2}_{\mathcal{T}_h}+\Delta t\int^{t^{n}}_{t^{n-1}}(\|\f{\p^{2} \phi}{\p t^{2}}\|^{2}_{\mathcal{T}_h}+\|\f{\p c_{f}}{\p t}\|^{2}_{\mathcal{T}_h})dt\}.
\end{equation}
\end{lem}
\begin{proof}

Denote the three terms on right-hand-side of \eqref{e21} by $I_{1}$,  $I_{2}$ and  $I_{3}$.   We estimate them one by one. 

For $I_{1}$, we have 
\[
I_{1}\leq C\|\bar{c}_{h}^{n-1}-{c}_{f}^{n-1}\|_{\mathcal{T}_h}\|q\|_{\mathcal{T}_h}\leq
C(\|\xi_{c}^{n-1}\|^{2}_{\mathcal{T}_h}+\|\eta_{c}^{n-1}\|^{2}_{\mathcal{T}_h}+\|q\|^{2}_{\mathcal{T}_h}).
\]
Notice that
\[
\begin{split}
I_{2}=&-(\f{\alpha K}{\rho_{s}}\phi^{n}_{h}(\bar{c}_{h}^{n-1}-{c}^{n-1}_{f}),q)_{\mathcal{T}_h}-(\f{\alpha K}{\rho_{s}}(\phi^{n}_{h}-\phi^{n}){c}^{n-1}_{f},q)_{\mathcal{T}_h}-(\f{\alpha K}{\rho_{s}}\phi^{n}({c}^{n-1}_{f}-{c}^{n}_{f}),q)_{\mathcal{T}_h}
\\
\leq &C\{\|\xi_{c}^{n-1}\|^{2}_{\mathcal{T}_h}+\|\eta_{c}^{n-1}\|^{2}_{\mathcal{T}_h}+\|q\|^{2}_{\mathcal{T}_h}+\Delta t\int^{t^{n}}_{t^{n-1}}\|\f{\p c_{f}}{\p t}\|^{2}_{\mathcal{T}_h}dt\}.
\end{split}
\]
Meanwhile, we know that
\[
I_{3}\leq C\{\Delta t\int^{t^{n}}_{t^{n-1}}(\|\f{\p^{2} \phi}{\p t^{2}}\|^{2}_{\mathcal{T}_h}+\|\f{\p c_{f}}{\p t}\|^{2}_{\mathcal{T}_h})dt+\|q\|^{2}_{\mathcal{T}_h}\}.
\]
Substituting the above estimates into \eqref{e21},  we get the inequality \eqref{e22}.
\end{proof}

\begin{lem}\label{lem9}
The approximate error of the discrete porosity satisfies
\begin{equation}\label{e23}
\|e_{\phi}^{n}\|^{2}_{\mathcal{T}_h}\leq C\{\Delta t\sum^{n}_{i=1}(\|\xi_{c}^{i-1}\|^{2}_{\mathcal{T}_h}+\|\eta_{c}^{i-1}\|^{2}_{\mathcal{T}_h})+\Delta t^{2}\int^{t^n}_{0}(\|\f{\p^{2} \phi}{\p t^{2}}\|^{2}_{\mathcal{T}_h}+\|\f{\p c_{f}}{\p t}\|^{2}_{\mathcal{T}_h})dt\}.
\end{equation}
\end{lem}
\begin{proof}
Taking $q=e^{n}_{\phi}$ in \eqref{e22}, we have
\begin{equation}\label{e24}
\begin{split}
&\f{1}{2\Delta t}[\|e_{\phi}^{n}\|^{2}_{\mathcal{T}_h}-\|e_{\phi}^{n-1}\|^{2}_{\mathcal{T}_h}]
\leq (\f{e_{\phi}^{n}-e_{\phi}^{n-1}}{\Delta t},e_{\phi}^{n})_{\mathcal{T}_h}
\\\leq&C\{\|\xi_{c}^{n-1}\|^{2}_{\mathcal{T}_h}+\|\eta_{c}^{n-1}\|^{2}_{\mathcal{T}_h}+\|e_{\phi}^{n}\|^{2}_{\mathcal{T}_h}+\Delta t\int^{t^{n}}_{t^{n-1}}(\|\f{\p^{2} \phi}{\p t^{2}}\|^{2}_{\mathcal{T}_h}+\|\f{\p c_{f}}{\p t}\|^{2}_{\mathcal{T}_h})dt\}.
\end{split}
\end{equation}
Multiplying \eqref{e23} by $2\Delta t$, for sufficiently small $\Delta t$, we have
\begin{equation}\label{e25}
\begin{split}
\|e_{\phi}^{n}\|^{2}_{\mathcal{T}_h}-\|e_{\phi}^{n-1}\|^{2}_{\mathcal{T}_h}\leq C\Delta t\{\|\xi_{c}^{n-1}\|^{2}_{\mathcal{T}_h}+\|\eta_{c}^{n-1}\|^{2}_{\mathcal{T}_h}+\|e_{\phi}^{n}\|^{2}_{\mathcal{T}_h}+\Delta t\int^{t^{n}}_{t^{n-1}}(\|\f{\p^{2} \phi}{\p t^{2}}\|^{2}_{\mathcal{T}_h}+\|\f{\p c_{f}}{\p t}\|^{2}_{\mathcal{T}_h})dt\}.
\end{split}
\end{equation}
Summing the above inequality \eqref{e25} over $n$, for sufficiently small $\Delta t$, we can get the estimate  \eqref{e23}.
\end{proof}


\subsection{The error estimates of pressure and velocity }

Using the definitions of projection operators and $B_{\mathbf{u}}$, we have
\[
B_{\mathbf{u}}(\phi^n_h; (\Pi^{RT}\mathbf{u}^n-\mathbf{u}^n,\Pi_h p^n-p^n,\Pi_e p^n-p^n),(\omega_h,v_h,\mu_h))
=(\beta(\phi^n_h)(\Pi^{RT}\mathbf{u}^n-\mathbf{u}^n),\omega_h)_{\mathcal{T}_h}.
\]
According to the boundedness and stability of  the bilinear form $B_{\mathbf{u}}$,  we have the estimate
\begin{equation}\nonumber
\begin{aligned}
&C_{\ast}\|(\Pi^{RT}\mathbf{u}^n-\mathbf{u}^n_h,\Pi_h p^n-p^n_h,\Pi_e p^n-\lambda^n_h)\|_\mathbf{u}\\
\leq&\sup\limits_{(\omega_h,v_h,\mu_h)}\frac{B_{\mathbf{u}}((\Pi^{RT}\mathbf{u}^n-\mathbf{u}^n_h,\Pi_h p^n-p^n_h,\Pi_e p^n-\lambda^n_h),(\omega_h,v_h,\mu_h))}{\| (\omega_h,v_h,\mu_h)\|_\mathbf{u}}\\
=&\sup\limits_{(\omega_h,v_h,\mu_h)}[\frac{B_{\mathbf{u}}((\Pi^{RT}\mathbf{u}^n-\mathbf{u}^n,\Pi_h p^n-p^n,\Pi_e p^n-\lambda^n),(\omega_h,v_h,\mu_h))}{\| (\omega_h,v_h,\mu_h)\|_\mathbf{u}}\\
&+\frac{(\f{\alpha K}{\rho_s}(1-\phi^n_h)(\bar{c}^{n-1}_h-c^n_f),v_h)_{\mathcal{T}_h}+(\f{\alpha K}{\rho_s}(\phi^n-\phi^n_h)c^n_f,v_h)_{\mathcal{T}_h}}{\| (\omega_h,v_h,\mu_h)\|_\mathbf{u}}]\\
=&\sup\limits_{(\omega_h,v_h,\mu_h)}[\frac{(\beta(\phi^n_h)(\Pi^{RT}\mathbf{u}^n-\mathbf{u}^n),\omega_h)_{\mathcal{T}_h}}{\| (\omega_h,v_h,\mu_h)\|_\mathbf{u}}\\
&+\frac{(\f{\alpha K}{\rho_s}(1-\phi^n_h)(\bar{c}^{n-1}_h-c^n_f),v_h)_{\mathcal{T}_h}+(\f{\alpha K}{\rho_s}(\phi^n-\phi^n_h)c^n_f,v_h)_{\mathcal{T}_h}}{\| (\omega_h,v_h,\mu_h)\|_\mathbf{u}}]\\
\leq&C(\| c^n_f-\bar{c}^{n-1}_h\|_{\mathcal{T}_h}+\| \Pi^{RT}\mathbf{u}^n-\mathbf{u}^n\|_{\mathcal{T}_h}+\| \phi^n-\phi^n_h\|_{\mathcal{T}_h}).
\end{aligned}
\end{equation}
Hence we get
\begin{equation}\label{e26}
\begin{aligned}
\|(\Pi^{RT}\mathbf{u}^n-\mathbf{u}^n_h,\Pi_h p^n-p^n_h,\Pi_e p^n-\lambda^n_h)\|_\mathbf{u}
\leq C\{\| c^n_f-\bar{c}^{n-1}_h\|_{\mathcal{T}_h}+\| \Pi^{RT}\mathbf{u}^n-\mathbf{u}^n\|_{\mathcal{T}_h}+\| \phi^n-\phi^n_h\|_{\mathcal{T}_h}\}.
\end{aligned}
\end{equation}
Using  \eqref{e26} and triangle inequality, we get the error estimate of the pressure and velocity.
\begin{lem}\label{lem10}
There exists  a constant $C$ independent of $\Delta t$ and $h$, such that the estimate 
\begin{equation}\label{e27}
\begin{aligned}
&\|(\mathbf{u}^n-\mathbf{u}^n_h,\Pi_h p^n-p^n_h,\Pi_e p^n-\lambda^n_h)\|_\mathbf{u}
\leq&C\{\| c^n_f-c^{n-1}_f\|_{\mathcal{T}_h}+\| c^{n-1}_f-c^{n-1}_h\|_{\mathcal{T}_h}+\|\Pi^{RT}\mathbf{u}^n-\mathbf{u}^n\|_{\mathcal{T}_h}+\| \phi^n-\phi^n_h\|_{\mathcal{T}_h}\}
\end{aligned}
\end{equation}
holds.
\end{lem}

\subsection{The error estimate of the concentration }
\begin{lem}\label{lem11}
Under the induction hypothesis \eqref{hyp1}, for any $n>0$, the following inequality holds:
\begin{equation}\label{e37}
\begin{aligned}
&\f{1}{\Delta t}[(\phi^n_h\xi_c^n,\xi_c^n)_{\mathcal{T}_h}-(\phi^{n-1}_h\xi_c^{n-1},\xi_c^{n-1})_{\mathcal{T}_h}] +\|\nabla\xi_c^n\|^2_{\mathcal{T}_h}+\f{1}{h}|\xi_c^n-\xi_\lambda^n|^2_{\partial\mathcal{T}_h}+2\langle|\mathbf{u}_h^{n}\cdot\mathbf{n}|(\xi_c^n-\xi_\lambda^n),\xi_c^n-\xi_\lambda^n\rangle_{\partial\mathcal{T}_h}
\\
\leq &C\{[1+h^{-d}\|\xi_\mathbf{u}^{n}\|^2_{\mathcal{T}_h}]\|\phi^n_h-\phi^n\|^2_{\mathcal{T}_h}+\|\xi_c^n\|^2_{\mathcal{T}_h}+\|\xi_\sigma^n\|^2_{\mathcal{T}_h}+\|\xi_{c}^{n-1}\|^{2}_{\mathcal{T}_h}+\|\mathbf{u}_h^{n}-\mathbf{u}^{n}\|^2_{\mathcal{T}_h}
\\
&+\|\eta_c^n\|^2_{\mathcal{T}_h}+\|\eta_c^{n-1}\|^2_{\mathcal{T}_h}+\|\eta_\sigma^n\|^2_{\mathcal{T}_h}+\frac{1}{\Delta t}\int^{t^n}_{t^{n-1}}\|\f{\p \eta_c}{\p t}\|^2_{\mathcal{T}_h}dt 
\\&+\Delta t\int^{t^{n}}_{t^{n-1}}(\|\f{\p^{2} \phi}{\p t^{2}}\|^{2}_{\mathcal{T}_h}+\|\f{\p c_{f}}{\p t}\|^{2}_{\mathcal{T}_h}+\|\f{\p^{2} }{\p t^{2}}(\phi c_f)\|^{2}_{\mathcal{T}_h})dt\}
\end{aligned}
\end{equation}
where $C$ is a constant independent of $h$ and $\Delta t$.
\end{lem}

\begin{proof}
For the convenience of analysis, we denote
\[
\begin{aligned}
&c^n_h-c^n_f=(c^n_h-\Pi_hc^n+(\Pi_hc^n_f-c^n_f)=\xi^n_c-\eta^n_c,\\
&\lambda^{c,n}_h-c^n_f=(\lambda^{c,n}_h-\Pi_ec^n_f)+(\Pi_ec^n_f-c^n_f)=\xi^n_\lambda-\eta^n_\lambda,\\
&\bm\sigma^n_h-\bm\sigma^n=(\bm\sigma^n_h-\Pi^{RT}\bm\sigma^n)+(\Pi^{RT}\bm\sigma^n-\bm\sigma^n)=\xi^n_{\sigma}-\eta^n_{\sigma}.
\end{aligned}
\]
We can derive the following formulas from Algorithm \ref{alg6}:
\begin{equation}\label{e28}
\begin{aligned}
&-\frac{1}{\Delta t}(\phi^n_h\xi_c^n,v_h)_{\mathcal{T}_h}-([K-f^n_P])\xi_c^n,v_h)_{\mathcal{T}_h}+(K\phi^n_h\xi_c^n,v_h)_{\mathcal{T}_h}
\\& +(\xi_\sigma^n,\nabla v_h)_{\mathcal{T}_h}+(\mathbf{u}_h^{n}\xi_c^n,\nabla v_h)_{\mathcal{T}_h}
+([\phi^n_hD(\mathbf{u}_h^{n})]^{-1}\xi_\sigma^n,\bm\tau_h)_{\mathcal{T}_h}\\
&+(\nabla\xi_c^n,\bm\tau_h)_{\mathcal{T}_h}+\langle\xi_\lambda^n,\bm\tau_h\cdot\mathbf{n}\rangle_{\partial\mathcal{T}_h}
-\langle\xi_c^n,\bm\tau_h\cdot\mathbf{n}\rangle_{\partial\mathcal{T}_h}+\langle\xi_\sigma^n\cdot\mathbf{n},\mu_h-v_h\rangle_{\partial\mathcal{T}_h}\\
&+\langle\mathbf{u}_h^{n}\cdot\mathbf{n}\xi_\lambda^n,\mu_h-v_h\rangle_{\partial\mathcal{T}_h^{in}}+\langle\mathbf{u}_h^{n}\cdot\mathbf{n}\xi_c^n,\mu_h-v_h\rangle_{\partial\mathcal{T}_h^{out}}\\
=&-\frac{1}{\Delta t}(\phi^{n-1}_h\xi_c^{n-1},v_h)_{\mathcal{T}_h}
-\frac{1}{\Delta t}(\phi^n_h\eta_c^n-\phi^{n-1}_h\eta_c^{n-1},v_h)_{\mathcal{T}_h}
\\&+\frac{1}{\Delta t}([\phi^n_h-\phi^n]c^n_f-[\phi^{n-1}_h-\phi^{n-1}]c^{n-1}_f,v_h)_{\mathcal{T}_h}+(\frac{\phi^nc^n_f-\phi^{n-1}c^{n-1}_f}{\Delta t}-\f{\p (\phi c_f)}{\p t},v_h)_{\mathcal{T}_h}
\\&-([K-f^n_P])\eta_c^n,v_h)_{\mathcal{T}_h}+(K\phi^n_h\eta_c^n,v_h)_{\mathcal{T}_h}-(K[\phi^n_h-\phi^n]c^n_f,v_h)_{\mathcal{T}_h}
\\&+(\mathbf{u}_h^{n}\eta_c^n,\nabla v_h)_{\mathcal{T}_h}-([\mathbf{u}_h^{n}-\mathbf{u}^{n}]c_f^n,\nabla v_h)_{\mathcal{T}_h}
\\&+([\phi^n_hD(\mathbf{u}_h^{n})]^{-1}\eta_\sigma^n,\bm\tau_h)_{\mathcal{T}_h}-(([\phi^n_hD(\mathbf{u}_h^{n})]^{-1}-[\phi^nD(\mathbf{u}^{n})]^{-1})\bm\sigma^n,\bm\tau_h)_{\mathcal{T}_h}
\\&+\langle\mathbf{u}_h^{n}\cdot\mathbf{n}\eta_c^n,\mu_h-v_h\rangle_{\partial\mathcal{T}_h}-\langle[\mathbf{u}_h^{n}-\mathbf{u}^{n}]\cdot\mathbf{n}c_f^n,\mu_h-v_h\rangle_{\partial\mathcal{T}_h}.
\end{aligned}
\end{equation}
Setting $v_h=-\xi_c^n$ and $\mu_h=-\xi_\lambda^n$,  $\bm\tau_h=\nabla\xi_c^n$ for every element $K\in\mathcal{T}_h$  and choosing $\bm\tau_h\cdot\mathbf{n}=\frac{1}{h}(\zeta_\lambda^n-\zeta_c^n)$ for $e\in\partial\mathcal{T}_h$ in \eqref{e28}, we can obtain
\begin{equation}\label{e29}
\begin{aligned}
&\frac{1}{\Delta t}(\phi^n_h\xi_c^n,\xi_c^n)_{\mathcal{T}_h}+([K-f^n_P])\xi_c^n,\xi_c^n)_{\mathcal{T}_h}-(K\phi^n_h\xi_c^n,\xi_c^n)_{\mathcal{T}_h}
\\& -(\xi_\sigma^n,\nabla \xi_c^n)_{\mathcal{T}_h}-(\mathbf{u}_h^{n}\xi_c^n,\nabla \xi_c^n)_{\mathcal{T}_h}
+([\phi^n_hD(\mathbf{u}_h^{n})]^{-1}\xi_\sigma^n,\nabla\xi_c^n)_{\mathcal{T}_h}\\
&+(\nabla\xi_c^n,\nabla\xi_c^n)_{\mathcal{T}_h}+\langle\xi_\lambda^n,\frac{1}{h}(\zeta_\lambda^n-\zeta_c^n)\rangle_{\partial\mathcal{T}_h}
-\langle\xi_c^n,\frac{1}{h}(\zeta_\lambda^n-\zeta_c^n)\rangle_{\partial\mathcal{T}_h}+\langle\xi_\sigma^n\cdot\mathbf{n},\xi_c^n-\xi_\lambda^n\rangle_{\partial\mathcal{T}_h}\\
&+\langle\mathbf{u}_h^{n}\cdot\mathbf{n}\xi_\lambda^n,\xi_c^n-\xi_\lambda^n\rangle_{\partial\mathcal{T}_h^{in}}+\langle\mathbf{u}_h^{n}\cdot\mathbf{n}\xi_c^n,\xi_c^n-\xi_\lambda^n\rangle_{\partial\mathcal{T}_h^{out}}\\
=&\frac{1}{\Delta t}(\phi^n_h\xi_c^n,\xi_c^n)_{\mathcal{T}_h}+(\nabla\xi_c^n,\nabla\xi_c^n)_{\mathcal{T}_h}+\f{1}{2}(\nabla\cdot\mathbf{u}_h^{n}\xi_c^n, \xi_c^n)_{\mathcal{T}_h}+([K-f^n_P])\xi_c^n,\xi_c^n)_{\mathcal{T}_h}-(K\phi^n_h\xi_c^n,\xi_c^n)_{\mathcal{T}_h}
\\&+\langle|\mathbf{u}_h^{n}\cdot\mathbf{n}|(\xi_c^n-\xi_\lambda^n),\xi_c^n-\xi_\lambda^n\rangle_{\partial\mathcal{T}_h} -(\xi_\sigma^n,\nabla \xi_c^n)_{\mathcal{T}_h}
+([\phi^n_hD(\mathbf{u}_h^{n})]^{-1}\xi_\sigma^n,\nabla\xi_c^n)_{\mathcal{T}_h}\\
&+\langle\xi_\lambda^n,\frac{1}{h}(\zeta_\lambda^n-\zeta_c^n)\rangle_{\partial\mathcal{T}_h}
-\langle\xi_c^n,\frac{1}{h}(\zeta_\lambda^n-\zeta_c^n)\rangle_{\partial\mathcal{T}_h}+\langle\xi_\sigma^n\cdot\mathbf{n},\xi_c^n-\xi_\lambda^n\rangle_{\partial\mathcal{T}_h}\\
=&
\frac{1}{\Delta t}(\phi^{n-1}_h\xi_c^{n-1},\xi_c^n)_{\mathcal{T}_h}
+\frac{1}{\Delta t}(\phi^n_h\eta_c^n-\phi^{n-1}_h\eta_c^{n-1},\xi_c^n)_{\mathcal{T}_h}
\\&+\frac{1}{\Delta t}([\phi^n_h-\phi^n]c^n_f-[\phi^{n-1}_h-\phi^{n-1}]c^{n-1}_f,\xi_c^n)_{\mathcal{T}_h}+(\frac{\phi^nc^n_f-\phi^{n-1}c^{n-1}_f}{\Delta t}-\f{\p (\phi c_f)}{\p t},\xi_c^n)_{\mathcal{T}_h}
\\&+([K-f^n_P])\eta_c^n,\xi_c^n)_{\mathcal{T}_h}-(K\phi^n_h\eta_c^n,\xi_c^n)_{\mathcal{T}_h}+(K[\phi^n_h-\phi^n]c^n_f,\xi_c^n)_{\mathcal{T}_h}
\\&-(\mathbf{u}_h^{n}\eta_c^n,\nabla \xi_c^n)_{\mathcal{T}_h}+([\mathbf{u}_h^{n}-\mathbf{u}^{n}]c_f^n,\nabla \xi_c^n)_{\mathcal{T}_h}
\\&+([\phi^n_hD(\mathbf{u}_h^{n})]^{-1}\eta_\sigma^n,\nabla\xi_c^n)_{\mathcal{T}_h}-(([\phi^n_hD(\mathbf{u}_h^{n})]^{-1}-[\phi^nD(\mathbf{u}^{n})]^{-1})\bm\sigma^n,\nabla\xi_c^n)_{\mathcal{T}_h}
\\&+\langle\mathbf{u}_h^{n}\cdot\mathbf{n}\eta_c^n,\xi_c^n-\xi_\lambda^n\rangle_{\partial\mathcal{T}_h}-\langle[\mathbf{u}_h^{n}-\mathbf{u}^{n}]\cdot\mathbf{n}c_f^n,\xi_c^n-\xi_\lambda^n\rangle_{\partial\mathcal{T}_h}.
\end{aligned}
\end{equation}
That is,
\begin{equation}\label{e30}
\begin{aligned}
&\frac{1}{\Delta t}(\phi^n_h\xi_c^n,\xi_c^n)_{\mathcal{T}_h}-\frac{1}{\Delta t}(\phi^{n-1}_h\xi_c^{n-1},\xi_c^n)_{\mathcal{T}_h}+(\nabla\xi_c^n,\nabla\xi_c^n)_{\mathcal{T}_h}
\\&+\langle\zeta_\lambda^n-\zeta_c^n,\frac{1}{h}(\zeta_\lambda^n-\zeta_c^n)\rangle_{\partial\mathcal{T}_h}+\langle|\mathbf{u}_h^{n}\cdot\mathbf{n}|(\xi_c^n-\xi_\lambda^n),\xi_c^n-\xi_\lambda^n\rangle_{\partial\mathcal{T}_h}
\\
=&\frac{1}{\Delta t}(\phi^n_h\eta_c^n-\phi^{n-1}_h\eta_c^{n-1},\xi_c^n)_{\mathcal{T}_h}
+\frac{1}{\Delta t}([\phi^n_h-\phi^n]c^n_f-[\phi^{n-1}_h-\phi^{n-1}]c^{n-1}_f,\xi_c^n)_{\mathcal{T}_h}
\\&-(\frac{(\phi c_f)^n-(\phi c_f)^{n-1}}{\Delta t}-\f{\p (\phi c_f)}{\p t},\xi_c^n)_{\mathcal{T}_h}
+(\xi_\sigma^n,\nabla \xi_c^n)_{\mathcal{T}_h}
-([\phi^n_hD(\mathbf{u}_h^{n})]^{-1}\xi_\sigma^n,\nabla\xi_c^n)_{\mathcal{T}_h}
\\&-([\f{1}{2}\nabla\cdot\mathbf{u}_h^{n}+K(1-\phi^n_h)-f^n_P]\xi_c^n, \xi_c^n)_{\mathcal{T}_h}
+([K-f^n_P])\eta_c^n,\xi_c^n)_{\mathcal{T}_h}-(K\phi^n_h\eta_c^n,\xi_c^n)_{\mathcal{T}_h}
\\&+(K[\phi^n_h-\phi^n]c^n_f,\xi_c^n)_{\mathcal{T}_h}
-(\mathbf{u}_h^{n}\eta_c^n,\nabla \xi_c^n)_{\mathcal{T}_h}+([\mathbf{u}_h^{n}-\mathbf{u}^{n}]c_f^n,\nabla \xi_c^n)_{\mathcal{T}_h}
\\&+([\phi^n_hD(\mathbf{u}_h^{n})]^{-1}\eta_\sigma^n,\nabla\xi_c^n)_{\mathcal{T}_h}-(([\phi^n_hD(\mathbf{u}_h^{n})]^{-1}-[\phi^nD(\mathbf{u}^{n})]^{-1})\bm\sigma^n,\nabla\xi_c^n)_{\mathcal{T}_h}
\\&-\langle\xi_\sigma^n\cdot\mathbf{n},\xi_c^n-\xi_\lambda^n\rangle_{\partial\mathcal{T}_h}
+\langle\mathbf{u}_h^{n}\cdot\mathbf{n}\eta_c^n,\xi_c^n-\xi_\lambda^n\rangle_{\partial\mathcal{T}_h}
\\&-\langle[\mathbf{u}_h^{n}-\mathbf{u}^{n}]\cdot\mathbf{n}c_f^n,\xi_c^n-\xi_\lambda^n\rangle_{\partial\mathcal{T}_h}
\\=&T_1+T_2+\cdots +T_{16}.
\end{aligned}
\end{equation}

It is easily seen that
\begin{equation}\label{e31}
\begin{aligned}
&\frac{1}{\Delta t}(\phi^n_h\xi_c^n,\xi_c^n)_{\mathcal{T}_h}-\frac{1}{\Delta t}(\phi^{n-1}_h\xi_c^{n-1},\xi_c^n)_{\mathcal{T}_h}
\\=&\f{1}{2\Delta t}[(\phi^n_h\xi_c^n,\xi_c^n)_{\mathcal{T}_h}-(\phi^{n-1}_h\xi_c^{n-1},\xi_c^{n-1})_{\mathcal{T}_h}]+\f{1}{2\Delta t}([\phi^n_h-\phi^{n-1}_h]\xi_c^n,\xi_c^n)_{\mathcal{T}_h}
\\&+\f{1}{2\Delta t}(\phi^{n-1}_h(\xi_c^n-\xi_c^{n-1}),\xi_c^n-\xi_c^{n-1})_{\mathcal{T}_h}
\\\geq&\f{1}{2\Delta t}[(\phi^n_h\xi_c^n,\xi_c^n)_{\mathcal{T}_h}-(\phi^{n-1}_h\xi_c^{n-1},\xi_c^{n-1})_{\mathcal{T}_h}]+\f{1}{2\Delta t}([\phi^n_h-\phi^{n-1}_h]\xi_c^n,\xi_c^n)_{\mathcal{T}_h}.
\end{aligned}
\end{equation}

Next, we will estimate the right-hand-side terms of \eqref{e30} one by one.   Using the boundedness of $\phi^n_h$ and $D$ and the regularity assumption of the solutiion, we know that
\begin{equation}\label{e32}
\begin{aligned}
T_4+\cdots+T_{13}\leq &C\{[1+h^{-d}\|\xi_\mathbf{u}^{n}\|^2_{\mathcal{T}_h}]\|\phi^n_h-\phi^n\|^2_{\mathcal{T}_h}+\|\xi_c^n\|^2_{\mathcal{T}_h}+\|\xi_\sigma^n\|^2_{\mathcal{T}_h}+\|\mathbf{u}_h^{n}-\mathbf{u}^{n}\|^2_{\mathcal{T}_h}\\
+&\|\eta_c^n\|^2_{\mathcal{T}_h}+\|\eta_\sigma^n\|^2_{\mathcal{T}_h}\}+\f{1}{2}\|\nabla\xi_c^n\|^2_{\mathcal{T}_h}
\end{aligned}
\end{equation}
where we have used  Young's inequality and the induction hypothesis
\begin{equation}\label{hyp1}
\|\nabla\cdot\mathbf{u}^n_h\|_{L^{\infty}(\mathcal{T}_h)}<C.
\end{equation}

Utilizing Schwarz inequality and the inverse inequality, we get the estimates of $T_{14}$, $T_{15}$ and $T_{16}$ as follows
\begin{equation}\label{e33}
\begin{aligned}
T_{14}+T_{15}+T_{16}\leq &C\{\|\xi_\sigma^n\|^2_{\mathcal{T}_h}+\|\mathbf{u}_h^{n}-\mathbf{u}^{n}\|^2_{\mathcal{T}_h}+\|\eta_c^n\|^2_{\mathcal{T}_h}\}+\f{1}{2h}|\xi_c^n-\xi_\lambda^n|^2_{\partial\mathcal{T}_h}.
\end{aligned}
\end{equation}

We  estimate $T_1$. By the proof of Theorem \ref{thm1}, we know that
\[
\f{\phi^n_h-\phi^{n-1}_h}{\Delta t}=\f{\alpha K (1-\phi^{n}_h)}{\rho_s}\bar{c}^{n-1}_h \leq C.
\]
So we get
\begin{equation}\label{e34}
\begin{aligned}
T_{1}= &\frac{1}{\Delta t}(\phi^n_h(\eta_c^n-\eta_c^{n-1}),\xi_c^n)_{\mathcal{T}_h} +\frac{1}{\Delta t}([\phi^n_h-\phi^{n-1}_h]\eta_c^{n-1},\xi_c^n)_{\mathcal{T}_h}\\
=&\frac{1}{\Delta t}(\phi^n_h\int^{t^n}_{t^{n-1}}\f{\p \eta_c}{\p t}dt,\xi_c^n)_{\mathcal{T}_h} +\frac{1}{\Delta t}([\phi^n_h-\phi^{n-1}_h]\eta_c^{n-1},\xi_c^n)_{\mathcal{T}_h}
\\\leq &\frac{1}{\Delta t}\int^{t^n}_{t^{n-1}}\|\f{\p \eta_c}{\p t}\|_{\mathcal{T}_h}dt\|\xi^n_c\|_{\mathcal{T}_h} +C\|\xi^n_c\|_{\mathcal{T}_h}\|\eta_c^{n-1}\|_{\mathcal{T}_h}
\\\leq &C(\frac{1}{\Delta t}\int^{t^n}_{t^{n-1}}\|\f{\p \eta_c}{\p t}\|^2_{\mathcal{T}_h}dt+\|\xi^n_c\|^2_{\mathcal{T}_h} +\|\eta_c^{n-1}\|^2_{\mathcal{T}_h}).
\end{aligned}
\end{equation}

For $T_2$, we have
\begin{equation}\label{e35}
\begin{aligned}
T_{2}= &\frac{1}{\Delta t}([\phi^n_h-\phi^n]c^{n-1}_f-[\phi^{n-1}_h-\phi^{n-1}]c^{n-1}_f+[\phi^n_h-\phi^n]\int^{t^n}_{t^{n-1}}\f{\p c_f}{\p t}dt,\xi_c^n)_{\mathcal{T}_h}
\\= &\frac{1}{\Delta t}([e_\phi^n-e_\phi^{n-1}]c^{n-1}_f,\xi_c^n)_{\mathcal{T}_h}+\frac{1}{\Delta t}([\phi^n_h-\phi^n]\int^{t^n}_{t^{n-1}}\f{\p c_f}{\p t}dt,\xi_c^n)_{\mathcal{T}_h}
\\\leq &C\{\|\xi_{c}^{n}\|^{2}_{\mathcal{T}_h}+\|\xi_{c}^{n-1}\|^{2}_{\mathcal{T}_h}+\|\eta_{c}^{n-1}\|^{2}_{{\mathcal{T}_h}}+\|\phi^n_h-\phi^n\|^{2}_{\mathcal{T}_h}
\\&+\Delta t\int^{t^{n}}_{t^{n-1}}(\|\f{\p^{2} \phi}{\p t^{2}}\|^{2}_{\mathcal{T}_h}+\|\f{\p c_{f}}{\p t}\|^{2}_{\mathcal{T}_h})dt\}.
\end{aligned}
\end{equation}
where we have used Lemma \ref{lem8} in the last inequality of \eqref{e35}.

For $T_3$, we can easily obtain
\begin{equation}\label{e36}
\begin{aligned}
T_{3}
\leq C\{\|\xi_{c}^{n}\|^{2}_{\mathcal{T}_h}+\Delta t\int^{t^{n}}_{t^{n-1}}\|\f{\p^{2} }{\p t^{2}}(\phi c_f)\|^{2}_{\mathcal{T}_h}dt\}.
\end{aligned}
\end{equation}
Substituting these estimates \eqref{e31}-\eqref{e36} into \eqref{e30}, we get the estimate of Lemma \ref{lem11}.

\end{proof}

From Lemma \ref{lem11}, the following estimate is easily obtained:

\begin{lem}\label{lem12}
Under Assumption 1, for any $n>0$, the following inequality holds:
\begin{equation}\label{e38}
\begin{aligned}
&\|\xi_c^n\|_{\mathcal{T}_h}^2+\Delta t\sum_{i=1}^n(\|\nabla\xi^i_c\|_{\mathcal{T}_h}^2+\f{1}{h}|\xi_c^i-\xi_\lambda^i|^2_{\partial\mathcal{T}_h})
\\
\leq &C\{\Delta t\sum^n_{i=1}([1+h^{-d}\|\xi_\mathbf{u}^{i}\|^2_{\mathcal{T}_h}]\|\phi^i_h-\phi^i\|^2_{\mathcal{T}_h}+\|\xi_c^i\|^2_{\mathcal{T}_h}+\|\xi_\sigma^i\|^2_{\mathcal{T}_h}+\|\xi_{c}^{i-1}\|^{2}_{\mathcal{T}_h}
\\&+\|\mathbf{u}_h^{i}-\mathbf{u}^{i}\|^2_{\mathcal{T}_h}
+\|\eta_c^i\|^2_{\mathcal{T}_h}+\|\eta_c^{i-1}\|^2_{\mathcal{T}_h}+\|\eta_\sigma^i\|^2_{\mathcal{T}_h})+\int^{t^n}_{0}\|\f{\p \eta_c}{\p t}\|^2_{\mathcal{T}_h}dt 
\\&+\Delta t^2\int^{t^{n}}_{0}(\|\f{\p^{2} \phi}{\p t^{2}}\|^{2}_{\mathcal{T}_h}+\|\f{\p c_{f}}{\p t}\|^{2}_{\mathcal{T}_h}+\|\f{\p^{2} }{\p t^{2}}(\phi c_f)\|^{2}_{\mathcal{T}_h})dt\}
\end{aligned}
\end{equation}
where $C$ is a constant independent of $h$ and $\Delta t$.
\end{lem}

\subsection{Convergence theorem}
\begin{thm}\label{thm3}
Suppose that  Assumptions \ref{ass1} holds, the coefficient $\beta(\phi)$ satisfies the Lipschitz continuousness with respect to $\phi$ and the solution of system \eqref{e1} has the regularity: $p\in L^2(0,T;H^{k+2}(\Omega))$,
$ \mathbf{u}\in L^\infty(0,T;H^{k+1}(\Omega))\cap L^\infty(0,T;W^{1,\infty}(\Omega))$, 
$c_f\in L^{2}(0,T;H^{k+1}(\Omega))\cap L^{\infty}(0,T;W^{1,\infty}(\Omega))$,
$\frac{\partial c_f}{\partial t}\in L^2(0,T;H^{k+1}(\Omega))\cap L^{\infty}(0,T;L^{\infty}(\Omega))$, $ \frac{\partial^2 \phi}{\partial t^2}, \frac{\partial^2 c_f}{\partial t^2}\in L^2(0,T;L^{2}(\Omega))$. And let $(\mathbf{u}^n_h, p^n_{h}, \lambda^n_{h},\bm\sigma^n_h,c^n_h,\lambda^{^n,c}_{h})$ be the solution of Algorithm \ref{alg6} with initial values $(\phi^0_h,c_h^{0})=(\Pi_h\phi_{0},\Pi_hc^{0}_f)$. Then, for some given $\Delta t_{0}>0$, when $\Delta t\leq \Delta t_{0}$, we have the following error estimate,  for $n>0$
\begin{equation}\label{e39}
\begin{aligned}
&\|c^n_h-c^n_f\|_{\mathcal{T}_h}^2+\|\phi^n_h-\phi^n\|_{\mathcal{T}_h}^2+\Delta t\sum^n_{i=1}(\|\bm\sigma_h^{i}-\bm\sigma^{i}\|^2_{\mathcal{T}_h}+\f{1}{h}|c_h^i-\lambda^i_h|^2_{\partial\mathcal{T}_h})
\leq C\{h^{2k+2}+\Delta t^2\},\\
&\|p^n_h-p^n\|_{\mathcal{T}_h}^2+\|\mathbf{u}^n_h-\mathbf{u}^n\|_{\mathcal{T}_h}^2\leq C\{h^{2k+2}+\Delta t^2\}.
\end{aligned}
\end{equation}
where $C$ is a constant independent of $h$ and $\Delta t$, and $1 \leq  k $ when $d = 2$, $2 \leq  k $ when $d = 3$.
\end{thm}
\begin{proof}
From \eqref{e11} and \eqref{e12}, we have the residual equation
\begin{equation}\label{e40}
\begin{aligned}
&(\up{a})\quad([\phi^n_hD(\mathbf{u}^{n}_h)]^{-1}\bm\sigma^{n}_h-[\phi^n D(\mathbf{u}^{n} )]^{-1}\bm\sigma^{n},\bm\tau_h)_{\mathcal{T}_h}-(c_h^{n}-c^n_f,\nabla\cdot\bm\tau_h)_{\mathcal{T}_h}+\langle\lambda^{c,n}_h,\bm\tau_h\cdot\mathbf{n}\rangle_{\partial\mathcal{T}_h}=0,\\
&(\up{b})\quad(\f{\phi^n_hc_h^{n}-\phi_h^{n-1}c_h^{n-1}-(\phi^n c_f^{n}-\phi^{n-1}c_f^{n-1})}{\Delta t},v_h)_{\mathcal{T}_h}
+(\nabla\cdot(\mathbf{u}^{n}_hc_h^{n})-\nabla\cdot(\mathbf{u}^{n}c_f^{n}),v_h)_{\mathcal{T}_h}
\\&\qquad\quad+(\nabla\cdot\bm\sigma_h^{n}-\nabla\cdot\bm\sigma^{n},v_h)_{\mathcal{T}_h}
+\langle\mathbf{u}_h^{n}\cdot\mathbf{n}(\lambda_h^{c,n}-c_h^{n})-\mathbf{u}^{n}\cdot\mathbf{n}(\lambda^{c,n}-c_f^{n}),v_h\rangle_{\partial\mathcal{T}_h^{in}}
\\&\qquad\quad+([K(1-\phi^n_h)-f^n_P]c^n_h-[K(1-\phi^n)-f^n_P]c^n_f,v_h)_{\mathcal{T}_h}
\\
&\qquad=(\f{\p (\phi c_f)}{\p t}-\f{\phi^n c_f^{n}-\phi^{n-1}c_f^{n-1}}{\Delta t},v_h)_{\mathcal{T}_h}.
\end{aligned}
\end{equation}
Choosing $(v_h,\bm\tau_h)=(\xi^n_c,\xi^n_\sigma)$  in \eqref{e40} and $\mu_h=\xi^n_\lambda$ in \eqref{e13}, we get 
\begin{equation}\label{e41}
\begin{aligned}
&\frac{1}{\Delta t}(\phi^n_h\xi_c^n,\xi_c^n)_{\mathcal{T}_h}-\frac{1}{\Delta t}(\phi^{n-1}_h\xi_c^{n-1},\xi_c^n)_{\mathcal{T}_h}
+([\phi^n_hD(\mathbf{u}_h^{n})]^{-1}\xi_\sigma^n,\xi_\sigma^n)_{\mathcal{T}_h}
\\
=&\frac{1}{\Delta t}(\phi^n_h\eta_c^n-\phi^{n-1}_h\eta_c^{n-1},\xi_c^n)_{\mathcal{T}_h}
-\frac{1}{\Delta t}([\phi^n_h-\phi^n]c^n_f-[\phi^{n-1}_h-\phi^{n-1}]c^{n-1}_f,\xi_c^n)_{\mathcal{T}_h}
\\&+(\f{\p (\phi c_f)}{\p t}-\f{\phi^n c_f^{n}-\phi^{n-1}c_f^{n-1}}{\Delta t},\xi_c^n)_{\mathcal{T}_h} -(\nabla\cdot(\mathbf{u}^{n}_hc_h^{n})-\nabla\cdot(\mathbf{u}^{n}c_f^{n}),\xi_c^n)_{\mathcal{T}_h}
\\&+([\phi^n_hD(\mathbf{u}_h^{n})]^{-1}\eta_\sigma^n,\xi_\sigma^n)_{\mathcal{T}_h}-(([\phi^n_hD(\mathbf{u}_h^{n})]^{-1}-[\phi^n D(\mathbf{u}^{n})]^{-1})\bm\sigma^n,\xi_\sigma^n)_{\mathcal{T}_h}
\\&
-([K(1-\phi^n_h)-f^n_P]c^n_h-[K(1-\phi^n)-f^n_P]c^n_f,\xi_c^n)_{\mathcal{T}_h}
\\&-\langle\mathbf{u}_h^{n}\cdot\mathbf{n}(\lambda_h^{c,n}-c_h^{n})-\mathbf{u}^{n}\cdot\mathbf{n}(\lambda^{c,n}-c_f^{n}),\xi_c^n\rangle_{\partial\mathcal{T}_h^{in}}.
\end{aligned}
\end{equation}
Taking $\mu_h=\xi^{n}_c$ in \eqref{e14}, we can get
\begin{equation}\nonumber
\begin{aligned}
\langle\mathbf{u}_h^{n}\cdot\mathbf{n}\hat{c}_h^{n},\xi^{n}_c\rangle_{\partial\mathcal{T}_h}=-\langle|\mathbf{u}_h^{n}\cdot\mathbf{n}|\lambda_h^{c,{n}},\xi^{n}_c\rangle_{\partial\mathcal{T}_h^{in}}+\langle|\mathbf{u}_h^{n}\cdot\mathbf{n}|c_h^{n},\xi^{n}_c\rangle_{\partial\mathcal{T}_h^{out}}=0.
\end{aligned}
\end{equation}
Thus we can reach
\begin{equation}\nonumber
\begin{aligned}
-\langle\mathbf{u}_h^{n}\cdot\mathbf{n}(\lambda_h^{c,{n}}-c_h^{n}),\xi^{n}_c\rangle_{\partial\mathcal{T}_h^{in}}
=&\langle|\mathbf{u}_h^{n}\cdot\mathbf{n}|(\lambda_h^{c,{n}}-c_h^{n}),\xi^{n}_c\rangle_{\partial\mathcal{T}_h^{in}}\\
=&\langle|\mathbf{u}_h^{n}\cdot\mathbf{n}|c_h^{n},\xi^{n}_c\rangle_{\partial\mathcal{T}_h^{out}}-\langle|\mathbf{u}_h^{n}\cdot\mathbf{n}|c_h^{n},\xi^{n}_c\rangle_{\partial\mathcal{T}_h^{in}}\\
=&\langle\mathbf{u}_h^{n}\cdot\mathbf{n}c_h^{n},\xi^{n}_c\rangle_{\partial\mathcal{T}_h^{out}}+\langle\mathbf{u}_h^{n}\cdot\mathbf{n}c_h^{n},\xi^{n}_c\rangle_{\partial\mathcal{T}_h^{in}}
=\langle\mathbf{u}_h^{n}\cdot\mathbf{n}c_h^{n},\xi^{n}_c\rangle_{\partial\mathcal{T}_h}.
\end{aligned}
\end{equation}
By use of Green's formula and the above equality, we know that
\begin{equation}\label{e42}
\begin{aligned}
&-(\nabla\cdot(\mathbf{u}^{n}_hc_h^{n})-\nabla\cdot(\mathbf{u}^{n}c_f^{n}),\xi_c^n)_{\mathcal{T}_h}-\langle\mathbf{u}_h^{n}\cdot\mathbf{n}(\lambda_h^{c,n}-c_h^{n})-\mathbf{u}^{n}\cdot\mathbf{n}(\lambda^{c,n}-c_f^{n}),\xi_c^n\rangle_{\partial\mathcal{T}_h^{in}}\\
=&(\nabla\cdot(\mathbf{u}^{n}c^{n})-\nabla\cdot(\mathbf{u}_h^{n}c_h^{n}),\xi^{n}_c)_{\mathcal{T}_h}-\langle\mathbf{u}^n\cdot\mathbf{n}-\mathbf{u}^n_h\cdot\mathbf{n},\xi^{n}_c\rangle_{\partial\mathcal{T}_h}\\
=&-(\mathbf{u}^{n}c^{n}-\mathbf{u}_h^{n}c_h^{n},\nabla\xi^{n}_c)_{\mathcal{T}_h}
=-(\mathbf{u}^{n}_{h}(\xi^{n}_c-\eta^{n}_c)+(\mathbf{u}^n-\mathbf{u}^n_h)c^{n},\nabla\xi^{n}_c)_{\mathcal{T}_h}\\
=&- (\mathbf{u}^{n}(\xi^{n}_c-\eta^{n}_c)+(\mathbf{u}^n-\mathbf{u}^n_h)c^{n},\nabla\xi^{n}_c)_{\mathcal{T}_h}
+ ([\mathbf{u}^{n}-\mathbf{u}^{n}_h](\xi^{n}_c-\eta^{n}_c)+(\mathbf{u}^n-\mathbf{u}^n_h)c^{n},\nabla\xi^{n}_c)_{\mathcal{T}_h}\\
\leq&C([1+h^{-d}\|\xi^{n}_\mathbf{u}\|_{\mathcal{T}_h}^2]\|\xi^n_c\|_{\mathcal{T}_h}^2+\|\eta^{n}_c\|_{\mathcal{T}_h}^2
+\|\mathbf{u}^n-\mathbf{u}^n_h\|_{\mathcal{T}_h}^2+\varepsilon\|\nabla\xi^{n}_c\|_{\mathcal{T}_h}^2),
\end{aligned}
\end{equation}
where we have used the fact that $\|\mathbf{u}^{n}\|_{L^{\infty}}, \|c^{n}\|_{L^{\infty}}\leq C$.

Substituting the estimate \eqref{e31}, \eqref{e32}, \eqref{e34}-\eqref{e36} and \eqref{e42} into \eqref{e40},  for sufficiently small $\varepsilon$, we can get the estimate with Lemma \ref{lem11}
\begin{equation}\label{e43}
\begin{aligned}
&\f{1}{\Delta t}[(\phi^n_h\xi_c^n,\xi_c^n)_{\mathcal{T}_h}-(\phi^{n-1}_h\xi_c^{n-1},\xi_c^{n-1})_{\mathcal{T}_h}] +\|\xi_\sigma^n\|^2_{\mathcal{T}_h}
\\
\leq &C\{[1+h^{-d}\|\xi_\mathbf{u}^{n}\|^2_{\mathcal{T}_h}](\|\phi^n_h-\phi^n\|^2_{\mathcal{T}_h}+\|\xi^n_c\|_{\mathcal{T}_h}^2)+\|\xi_{c}^{n-1}\|^{2}_{\mathcal{T}_h}
\\
&+\|\mathbf{u}_h^{n}-\mathbf{u}^{n}\|^2_{\mathcal{T}_h}
+\|\eta_c^n\|^2_{\mathcal{T}_h}+\|\eta_c^{n-1}\|^2_{\mathcal{T}_h}+\|\eta_\sigma^n\|^2_{\mathcal{T}_h}+\frac{1}{\Delta t}\int^{t^n}_{t^{n-1}}\|\f{\p \eta_c}{\p t}\|^2_{\mathcal{T}_h}dt 
\\&+\Delta t\int^{t^{n}}_{t^{n-1}}(\|\f{\p^{2} \phi}{\p t^{2}}\|^{2}_{\mathcal{T}_h}+\|\f{\p c_{f}}{\p t}\|^{2}_{\mathcal{T}_h}+\|\f{\p^{2} }{\p t^{2}}(\phi c_f)\|^{2}_{\mathcal{T}_h})dt\} .
\end{aligned}
\end{equation}
Multiplying the above estimate by $\Delta t$ and summing it over $n$, we can get the estimate with Lemma \ref{lem9} and Lemma \ref{lem10} as follows
\begin{equation}\label{e44}
\begin{aligned}
&\|\xi_c^n\|^2_{\mathcal{T}_h} +\Delta t\sum^n_{i=1}\|\xi_\sigma^i\|^2_{\mathcal{T}_h}
\\
\leq &C\{\Delta t\sum^n_{i=1}([1+h^{-d}\max_{0\leq j\leq i}\|\xi_c^{i-1}\|^2_{\mathcal{T}_h}]\|\xi^i_c\|_{\mathcal{T}_h}^2)+\|\xi_{c}^{i-1}\|^{2}_{\mathcal{T}_h}
+\|\eta_c^i\|^2_{\mathcal{T}_h}+\|\eta_c^{i-1}\|^2_{\mathcal{T}_h}+\|\eta_\sigma^i\|^2_{\mathcal{T}_h})
\\&+\int^{t^n}_{0}\|\f{\p \eta_c}{\p t}\|^2_{\mathcal{T}_h}dt 
+\Delta t^2\int^{t^{n}}_{0}(\|\f{\p^{2} \phi}{\p t^{2}}\|^{2}_{\mathcal{T}_h}+\|\f{\p c_{f}}{\p t}\|^{2}_{\mathcal{T}_h}+\|\f{\p^{2} }{\p t^{2}}(\phi c_f)\|^{2}_{\mathcal{T}_h})dt\}.
\end{aligned}
\end{equation}

In order to complete our proof, we need the following inductive hypothesis
\begin{equation}\label{e45}
 h^{-d/2}\max_{0\leq i\leq n-1}\|\xi^i_c\|_{\mathcal{T}_h} \leq C.
\end{equation}
Using Gronwall's lemma, we have the estimate
\begin{equation}\label{e46}
\|\xi_c^n\|^2_{\mathcal{T}_h} +\Delta t\sum^n_{i=1}\|\xi_\sigma^i\|^2_{\mathcal{T}_h}
\leq C\{h^{2s+2}+\Delta t^2\}.
\end{equation}
Using Lemmas \ref{lem9}, \ref{lem10}, \ref{lem11} and  \ref{lem12}, we  get the first estimate of \eqref{e40} and
\begin{equation}\label{e47}
\begin{aligned}
\|\mathbf{u}^n-\mathbf{u}^n_h\|^2_{\mathcal{T}_h} + \|\nabla(\Pi_hp^n-p^n_h)\|^2_{\mathcal{T}_h}
\leq C\{h^{2s+2}+\Delta t^2\}.
\end{aligned}
\end{equation}

As we know, the above error estimates  are obtained under the inductive hypotheses \eqref{hyp1} and \eqref{e45}.
Now we check it.   from \eqref{e47} we know that
\begin{equation}\nonumber
\begin{aligned}
\|\nabla\cdot \mathbf{u}^n_h\|_{L^\infty} \leq 
&\|\nabla\cdot (\mathbf{u}^n_h-\Pi^{RT}\mathbf{u}^n)\|_{L^\infty(\mathcal{T}_h)} +\|\Pi^{RT}\mathbf{u}^n\|_{L^\infty(\mathcal{T}_h)}
\\ \leq &Ch^{-\f{d}{2}}\|\nabla\cdot (\mathbf{u}^n_h-\Pi^{RT}\mathbf{u}^n)\|_{L^2(\mathcal{T}_h)} +\|\Pi^{RT}\mathbf{u}^n\|_{L^\infty(\mathcal{T}_h)}
 \\
 \leq &Ch^{-\f{d}{2}-1}\|\mathbf{u}^n_h-\Pi^{RT}\mathbf{u}^n\|_{L^2(\mathcal{T}_h)} +\|\Pi^{RT}\mathbf{u}^n\|_{L^\infty(\mathcal{T}_h)}
 \\
 \leq &Ch^{-\f{d}{2}-1} (h^{s+1}+\Delta t)+\|\Pi^{RT}\mathbf{u}^n\|_{L^\infty(\mathcal{T}_h)},
\end{aligned}
\end{equation}
and
\begin{equation}\nonumber
h^{-\f{d}{2} }\|\xi^{n}_c\|_{\mathcal{T}_h}\leq Ch^{-d/2}(h^{s+1}+\Delta t)\leq C.
\end{equation}
Hence the induction hypotheses holds.

Finally, we estimate the boundedness of $\| p^n-p^n_h\|_{\mathcal{T}_h}$. As in \cite{BF1991}, we know that, if $(p^n_h,\mathbf{u}^n_h, \lambda^n_h)$ is the solution of  the hybrid
mixed finite element method \eqref{e8},  then $(p_h,\mathbf{u}_h)$ is the solution of the classical mixed finite element method:
\[
\begin{aligned}
\left(\frac{\alpha K}{\rho_s}(1-\phi^n_h)\bar{c}^{n-1}_h,v_h\right)_{\mathcal{T}_h}+(\nabla\cdot \mathbf{u}^n_h,v_h)_{\mathcal{T}_h}&=(f^n,v_h)_{\mathcal{T}_h},\\
\left(\beta(\phi^n_h)\mathbf{u}^n_h,\bm\tau_h\right)_{\mathcal{T}_h}-(p^n_h,\nabla\cdot\bm\tau_h)_{\mathcal{T}_h}&=0.
\end{aligned}
\]
The mixed variational weak formulation of the velocity and pressure can be written as follows
\[
\begin{aligned}
\left(\frac{\alpha K}{\rho_s}(1-\phi^n){c^n}_f,v\right)_{\mathcal{T}_h}+(\nabla\cdot \mathbf{u}^n,v)_{\mathcal{T}_h}&=(f,v)_{\mathcal{T}_h},\\
\left(\beta(\phi^n)\mathbf{u}^n,\bm\tau\right)_{\mathcal{T}_h}-(p^n,\nabla\cdot\bm\tau)_{\mathcal{T}_h}&=0.
\end{aligned}
\]
Hence, we have
\[
\begin{aligned}
(p^n_h-\Pi_h p^n,\nabla\cdot\omega_h)_{\mathcal{T}_h}=&(p^n-\Pi_h p^n,\nabla\cdot\bm\tau_h)_{\mathcal{T}_h}+(\beta(\phi^n_h)\mathbf{u}^n_h-\beta(\phi^n)\mathbf{u}^n,\bm\tau_h)_{\mathcal{T}_h}
\\=&(\beta(\phi^n_h)(\mathbf{u}^n_h-\mathbf{u}^n),\bm\tau_h)_{\mathcal{T}_h}+([\beta(\phi^n_h)-\beta(\phi^n)]\mathbf{u}^n,\bm\tau_h)_{\mathcal{T}_h}
\end{aligned}
\]
where we have used the definition of the projection operator $\Pi_h$. By the inf-sup condition of the mixed finite element spaces, we get
\[
\begin{aligned}
C_0 \|p^n_h-\Pi_h p^n\|_{\mathcal{T}_h}\leq & 
\sup_{\bm\tau_h \in W_h} \f{(p^n_h-\Pi_h p^n,\nabla\cdot\bm\tau_h)_{\mathcal{T}_h}}{\|\omega_h\|_{\Theta_h}}
\leq C\{\|\mathbf{u}^n_h-\mathbf{u}^n\|_{\mathcal{T}_h}+\|\phi^n_h-\phi^n\|_{\mathcal{T}_h}\}
\end{aligned}
\]
where $C_0>0$ is a constant independent of $h$. Using  \eqref{e47}, we complete our proof of Theorem \ref{thm3}.
\end{proof}
\section{Numerical examples}
\setcounter{equation}{0}
\subsection{Convergence test}
\subsubsection{Elliptic type pressure equation}
In this study, we examine equation \eqref{e2}(a,b) in conjunction with Algorithm \ref{alg3}. The domain is defined as $\Omega=[0,1]\times[0,1]$, with parameters set to $K=0$ and $\kappa=\mu=1$. The exact pressure function is specified as $p=\sin(\pi x)\sin(\pi y)$, from which the source function $f$ and boundary conditions are derived. For varying mesh sizes $h$ and polynomial degrees $k$ of the Raviart-Thomas mixed finite element space, numerical results for both $p$ and $\mathbf{u}$ are presented in Tables \ref{tab1} and \ref{tab2}, respectively. These results demonstrate that Algorithm \ref{alg3} achieves optimal convergence accuracy.

\renewcommand\arraystretch{1.25}
\begin{table}[!ht]
\centering 
\caption{ Numerical results  of  the pressure $p$ in $L^2 $-norm for  elliptic type pressure equation}\label{tab1}
\begin{tabular}{c|cccccc}
\bottomrule
$h$&$k=0$&\text{Rate}&$k=1$&\text{Rate}&$k=2$&\text{Rate}\\
\hline
0.1&$4.6542e{-02}$&-&$2.5513e{-03}$&-&$7.9528e{-05}$&- \\
0.05&$2.2510e{-02}$& 1.0480  &$5.6873e{-04}$& 2.1654 &$8.5868e{-06}$& 3.2113\\
0.025&$1.1134e{-02}$& 1.0156  &$1.3741e{-04}$&2.0493  &$1.0838e{-06}$&2.9860 \\
\bottomrule
\end{tabular}
\end{table}
\begin{table}[!ht]
\centering 
\caption{Numerical results  of the Darcy velocity $\bf{u}$ in $ L^2 $-norm for  elliptic type pressure equation}\label{tab2}
\begin{tabular}{ccccccc}
\bottomrule
$h$&$k=0$&\text{Rate}&$k=1$&\text{Rate}&$k=2$&\text{Rate}\\
\hline
0.1&$2.0567e{-01} $&-&$7.4518e{-03} $ &-&$2.6019e{-04}$ &- \\
0.05&$1.0093e{-01} $& 1.0270 &$1.7083e{-03}$ &2.1250    &$2.7213e{-05} $&3.2572  \\
0.025&$5.0553e{-02}$&0.9975   &$4.2955e{-04}$ & 1.9917  &$3.2133e{-06} $& 3.0822 \\
\bottomrule
\end{tabular}
\end{table}

\subsubsection{Convection-diffusion type concentration equation}
Here we consider the convection-diffusion problem \eqref{e2}(c) using Algorithm \ref{alg5}. To simplify computations,  we set $\mathbf{u}=[y,x]^T$, $\phi=1$, $K=0$, $f_P=-x$ and $c_I=1$. The exact concentration  is defined as $c_f=e^{t-x-y}\sin(\pi x) \sin(\pi y)$,  and the function $f_I$ along with initial-boundary conditions are derived from the exact solution. For varying diffusion coefficient $D$, mesh size $h$ and polynomial degrees $k$, we compute the numerical errors and convergence rates, which are summarized in in Table \ref{tab3} and Table \ref{tab4}, where the time increment  is set to  $\Delta t=1e-4$. These results indicate that Algorithm \ref{alg5}  exhibits optimal convergence properties.
 \renewcommand\arraystretch{1.25}
\begin{table}[!ht]
\centering 
\caption{Numerical results  of the concentration $c_f$ and the flux $\bm\sigma$ in $ L^2 $-norm for  concentration equation with $D=1.0$}\label{tab3}
\begin{tabular}{c|cccc|cccc}
\bottomrule
\multirow{2}{*}{$h$}&
\multicolumn{4}{c|}{$\|c_h-c_f\|_{L^2}$}&\multicolumn{4}{c}{$\|\bm\sigma-\bm\sigma_h\|_{L^2}$}\\
\cline{2-9}
&$k=0$&\text{Rate}&$k=1$&\text{Rate}&$k=0$&\text{Rate}&$k=1$&\text{Rate}\\
\midrule
$0.1$&$5.1558e{-02}$ &-&$3.2333e{-03} $&- &$2.6777e{-01} $&-&$1.3469e{-02} $ &-\\
$0.05$&$2.5285e{-02}$  & 1.0279 &$7.1038e{-04} $&2.1864   &$1.3206e{-01}$& 1.0198&$3.2289e{-03}$  &2.0605\\
$0.025$&$1.2653{e-02}$ &  0.9988&$1.7401e{-04}$& 2.0294 &$6.5342e{-02}$& 1.0151 &$7.8314e{-04}$ &2.0437 \\
\bottomrule
\end{tabular}
\end{table}

 \renewcommand\arraystretch{1.25}
\begin{table}[!ht]
\centering 
\caption{Numerical results  of the concentration $c_f$ and the flux $\bm\sigma$ in $ L^2 $-norm for  concentration equation with $D=0.01$}\label{tab4}
\begin{tabular}{c|cccc|cccc}
\bottomrule
\multirow{2}{*}{$h$}&
\multicolumn{4}{c|}{$\|c_h-c_f\|_{L^2}$}&\multicolumn{4}{c}{$\|\bm\sigma-\bm\sigma_h\|_{L^2}$}\\
\cline{2-9}
&$k=0$&\text{Rate}&$k=1$&\text{Rate}&$k=0$&\text{Rate}&$k=1$&\text{Rate}\\
\midrule
$0.1$&$6.0094e-02$&-              & $3.3389e-03$&-                 &$9.4209e-03	$&-&$7.8256e-04$&-      \\
$0.05$&$2.9210e-02$&$1.0407 $&$7.2005e-04$&$2.2132 $&$5.9263e-03$&0.6687&$2.2495e-04$&$1.7986 $ \\
$0.025$&$1.4549e-02$&$1.0055 $&$1.8901e-04$&$1.9296 $&$3.3895e-03$&0.8061&$6.2568e-05$&$1.8461 $ \\
\bottomrule
\end{tabular}
\end{table}

\subsubsection{Coupled problem}
In this study, we maintain the domain as  $\Omega=[0,1]\times[0,1]$, and select the parameters accordingly:
\[
T=1.0,\,\up{D}=d\mathbf{I},\, d=1.0e-2, \, \kappa_0=1,\, \alpha=\kappa_c=\kappa_s=\mu=f_I=1,\,
a_0=0.5,\, \rho_s=10,
\]
where $\mathbf{I}$ is an identity matrix. The initial-boundary-value conditions and the right-hand-side are determined based on the exact solutions.
\[
\left\{
\begin{split}
&p(x,y,t)=t\sin(\pi x)\sin(\pi y),
\\
&c_f(x,y,t)= tx^2(1-x)^2y^2(1-y)^2,\\
&\phi(x,y,t)=1-\exp\{{-\f{1}{80}t^2x^2(1-x)^2y^2(1-y)^2e^{(x+y+1)}-(x+y+1)}\}, \\
\end{split}
\right.
\]
To verify the convergence accuracy of our proposed method, we conduct numerical experiments with a fixed time step $\Delta t=1e-4$ and varying spatial steps $h$. The results are summarized in Table \ref{tab5}. From these results, it is evident that our method achieves second-order accuracy when $k=1$, which aligns with our theoretical analysis.

\begin{table}
\caption{Numerical results for coupled problem. }\label{tab5}
\begin{center}
\begin{tabular}{c|cccccccc}
\bottomrule
$h$& $\|p-p_h\|_{L^2}$&Rate &$\|\mathbf{u}-\mathbf{u}_h\|_{L^2}$&Rate  &  $\|\phi-\phi_{h}\|_{L^{2}}$ & Rate &$\|c_{f}-c_{h}\|_{L^{2}}$ & Rate 
\\
\hline
$0.2$        &$2.0301e-02$&-              &$1.7355e-01$&-&$1.7380e-03$&-&$1.9424e-04$&-\\
$0.1$        &$5.0597e-03$&$2.0044$ &$4.5963e-02$& $1.9168$&$4.4233e-04$&$2.0342$&$5.3631e-05$&$1.8567$
\\
$0.05$        &$1.2981e-03$& $1.9626$ & $1.1166e-02$&$1.9792$&$9.0000e-05$& $2.0997$&$1.3846e-05$&$1.9536$\\
$0.025$        &$3.2826e-04$&$1.9835$ &$2.9251e-03$&$1.9947$&$2.3221e-05$&  $1.9835$&$3.5442e-06$&$1.9659$
 \\
\bottomrule
\end{tabular}
\end{center}
\end{table}

\bigskip

\subsection{Simulation of wormhole propagation}

Here we simulate the propagation of wormholes within a rectangular rock tube. The domain is defined as $\Omega=[0,0.2m]\times[0,0.2m]$ and the physical parameters are set as follows
\[
\begin{split}
& \mu=0.01Pas,\quad \rho_s=2500kg/m^2,\quad \alpha=100kg/mol,\quad a_0=2m^{-1},\quad \kappa_0=10^{-9}m^2,
\\&
\up{D}=d\mathbf{I},\quad d=1e-5,\quad\kappa_c=1m/s,\quad \kappa_s=10m/s,\quad \phi_0=0.3,\quad c_I=1mol/m^2, \quad c^0_f=0.
\end{split}
\]
To investigate the wormhole propagation, the initial porosity and permeability are defined as
\[
\left\{\begin{split}
&\phi_0=0.6,\quad (x,y)=(0, 0.1),\\
&\phi_0=0.3, \quad\textrm{otherwise},
\end{split}
\right.\quad
\left\{\begin{split}
&\kappa_0=10^{-7},\quad (x,y)=(0, 0.1),\\
&\kappa_0=10^{-9}, \quad\textrm{otherwise}.
\end{split}
\right.
\]
The boundary conditions are specified such that the top and bottom boundaries satisfy $\mathbf{u}\cdot\mathbf{n}=0$. Acid is injected into the porous medium from the left boundary at a velocity of $0.02m/s$ and is drained from the right boundary at the same velocity.   The injection and production flow rates are given by:
\[
f_I=\left\{
\begin{split}
1,&\quad (x,y)=(0,0.1)\\
0,&\quad\textrm{otherwise}
\end{split}\right.
,\quad
f_P=\left\{
\begin{split}
-1,&\quad (x,y)=(0.2, 0.1)\\
0,&\quad\textrm{otherwise}
\end{split}\right.
\]

\begin{figure}
\centering
\subfigure[]{
    \label{fig1:a}
\includegraphics[width=6.5cm,height=5cm]{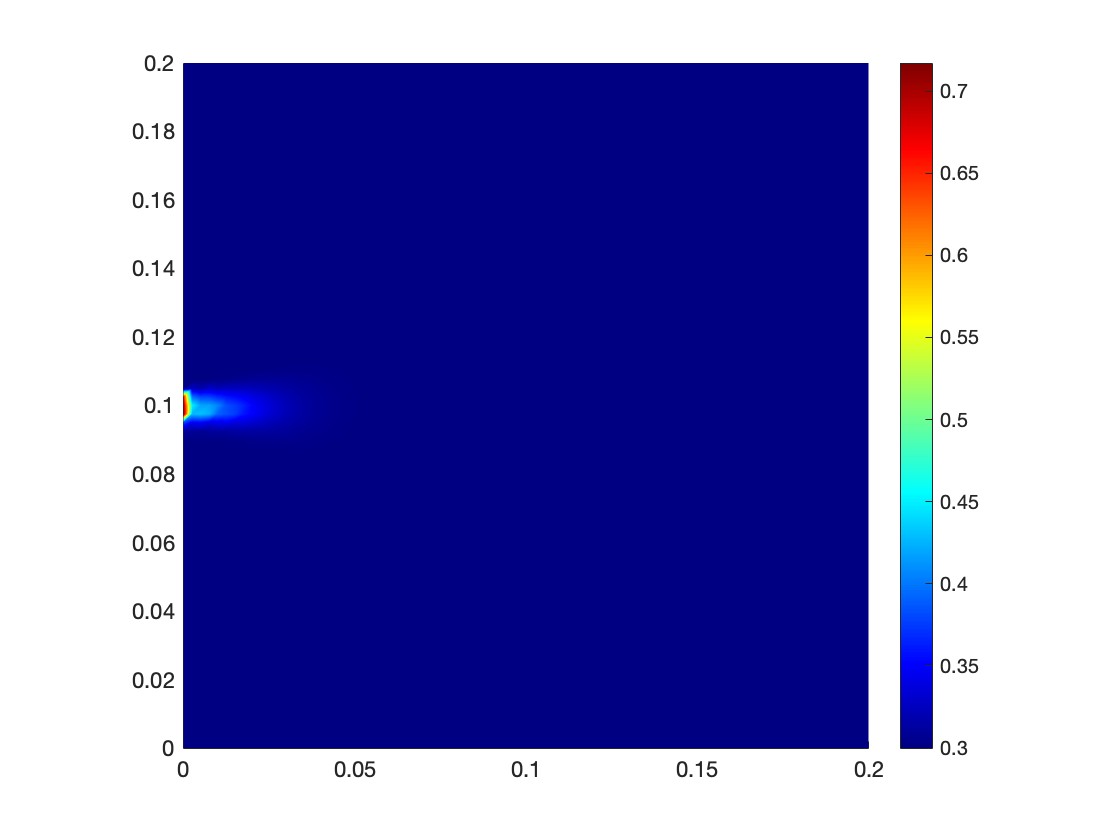}}
\subfigure[]{
    \label{fig1:b}
\includegraphics[width=6.5cm,height=5cm]{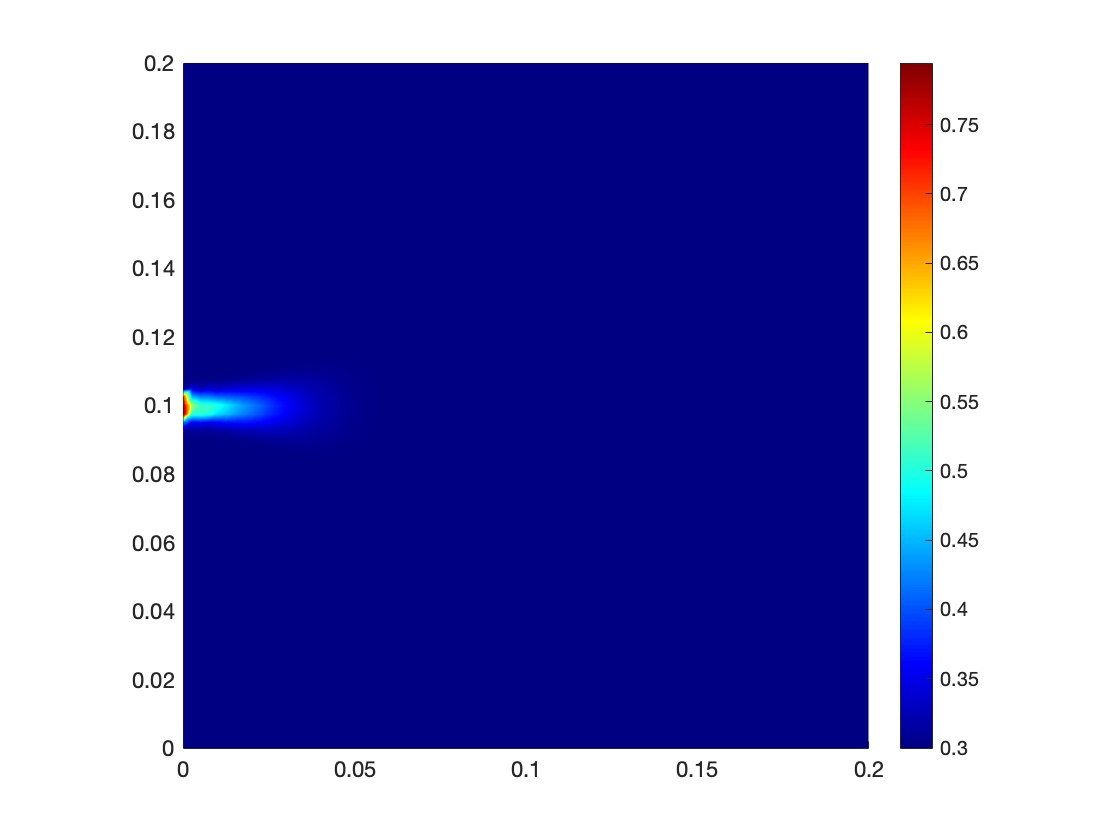}}
\subfigure[]{
    \label{fig1:c}
\includegraphics[width=6.5cm,height=5cm]{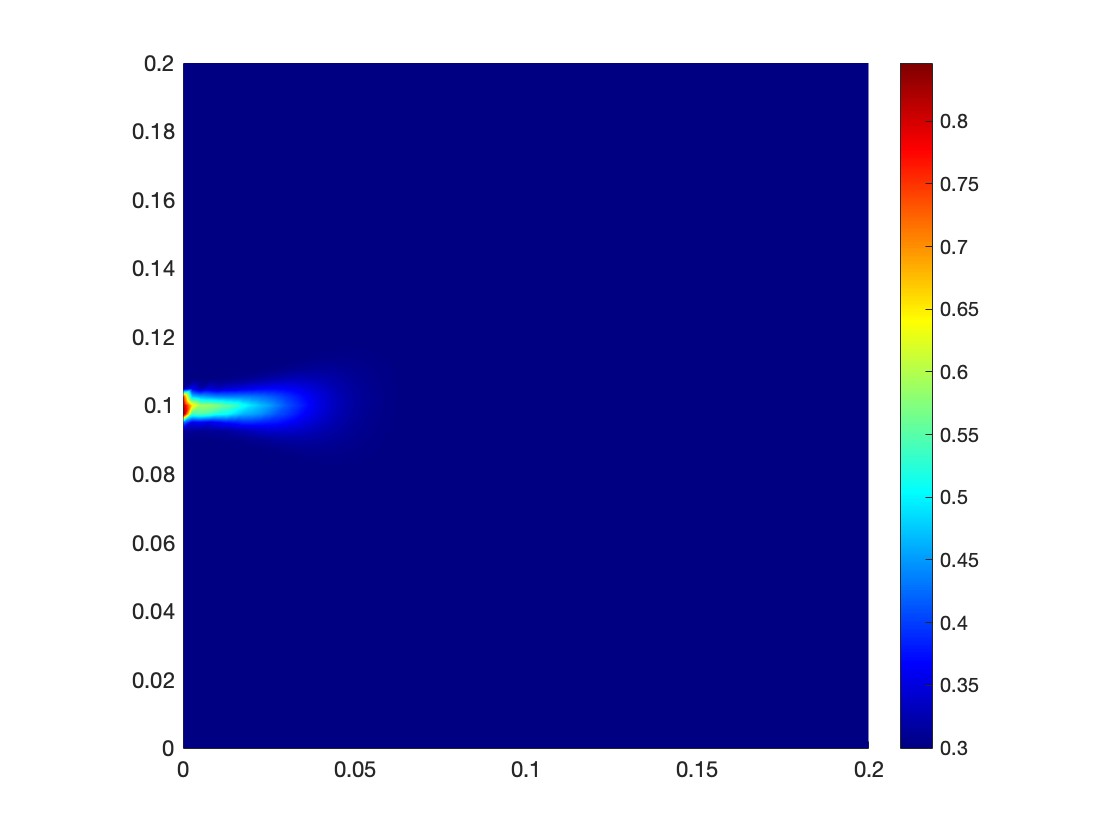}}
\subfigure[]{
    \label{fig1:d}
\includegraphics[width=6.5cm,height=5cm]{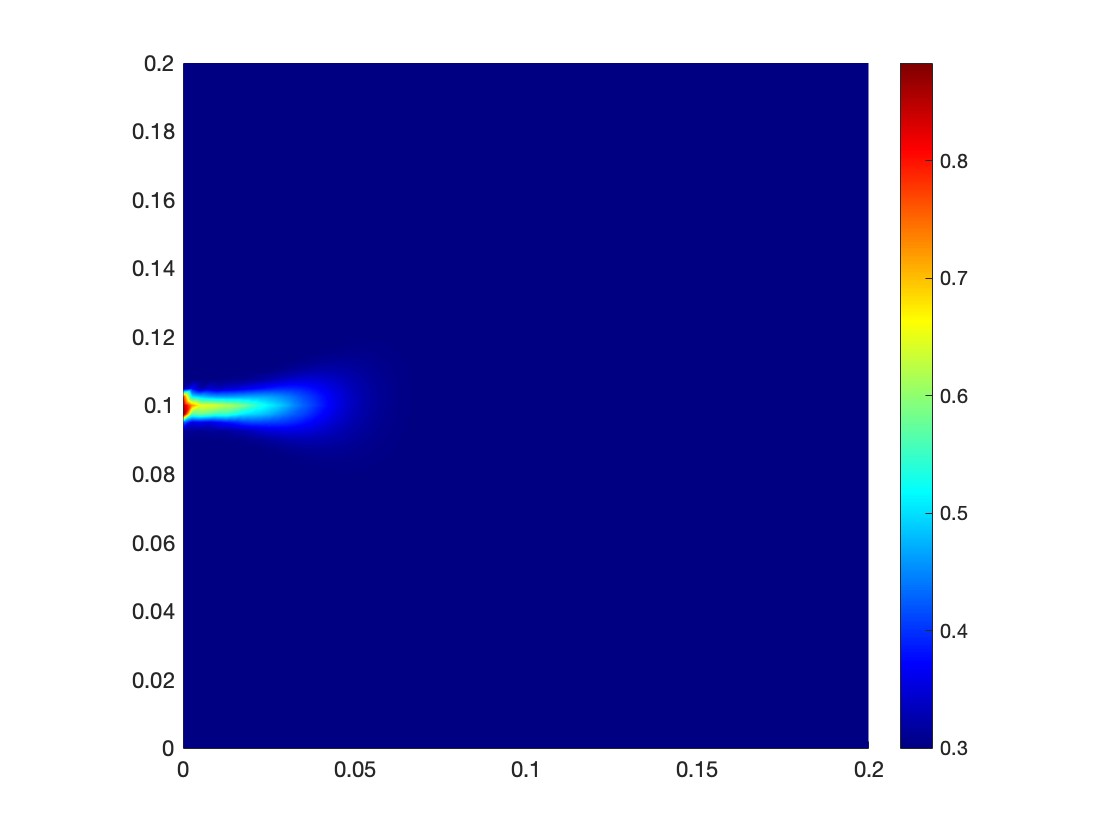}}
\caption{Porosity of rock at time $t=10s$ (a), $t=20s$ (b), $t=30s$ (c),
$t=40s$ (d) with $d=1e-5$.}
\label{fig1}
\end{figure}

\begin{figure}
\centering
\subfigure[]{
    \label{fig2:a}
\includegraphics[width=6.5cm,height=5cm]{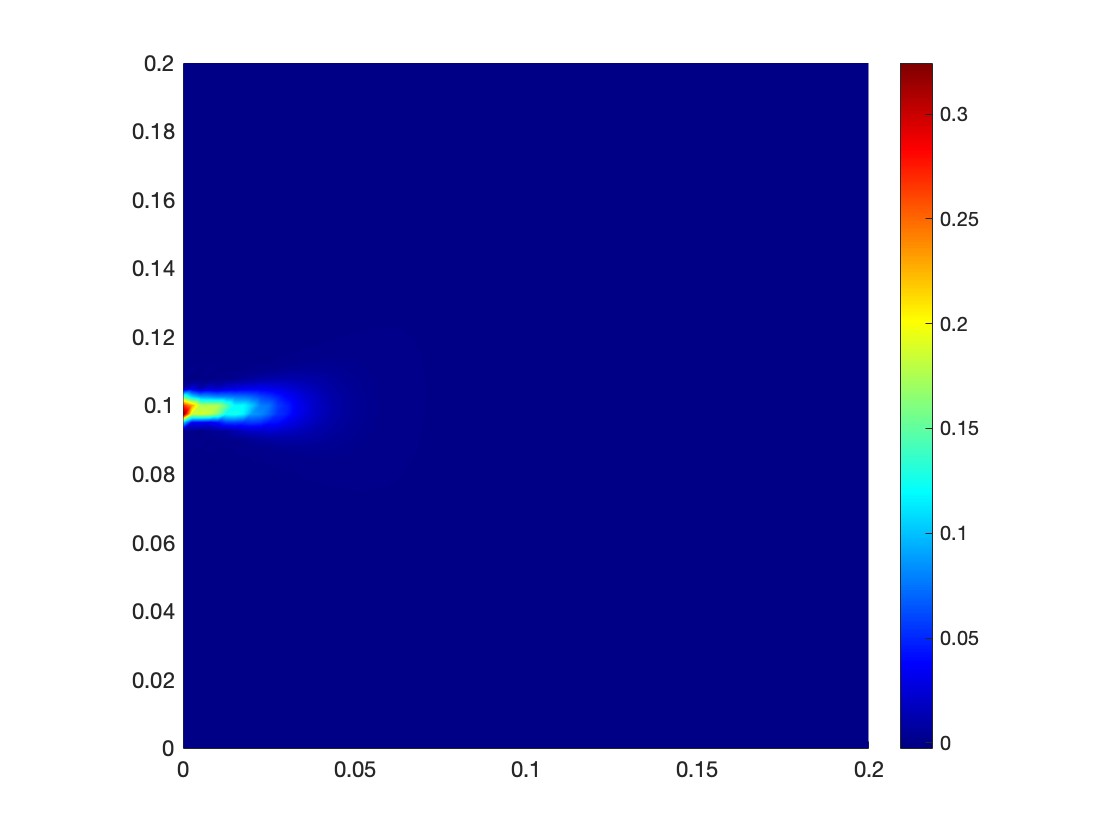}}
\subfigure[]{
    \label{fig2:b}
\includegraphics[width=6.5cm,height=5cm]{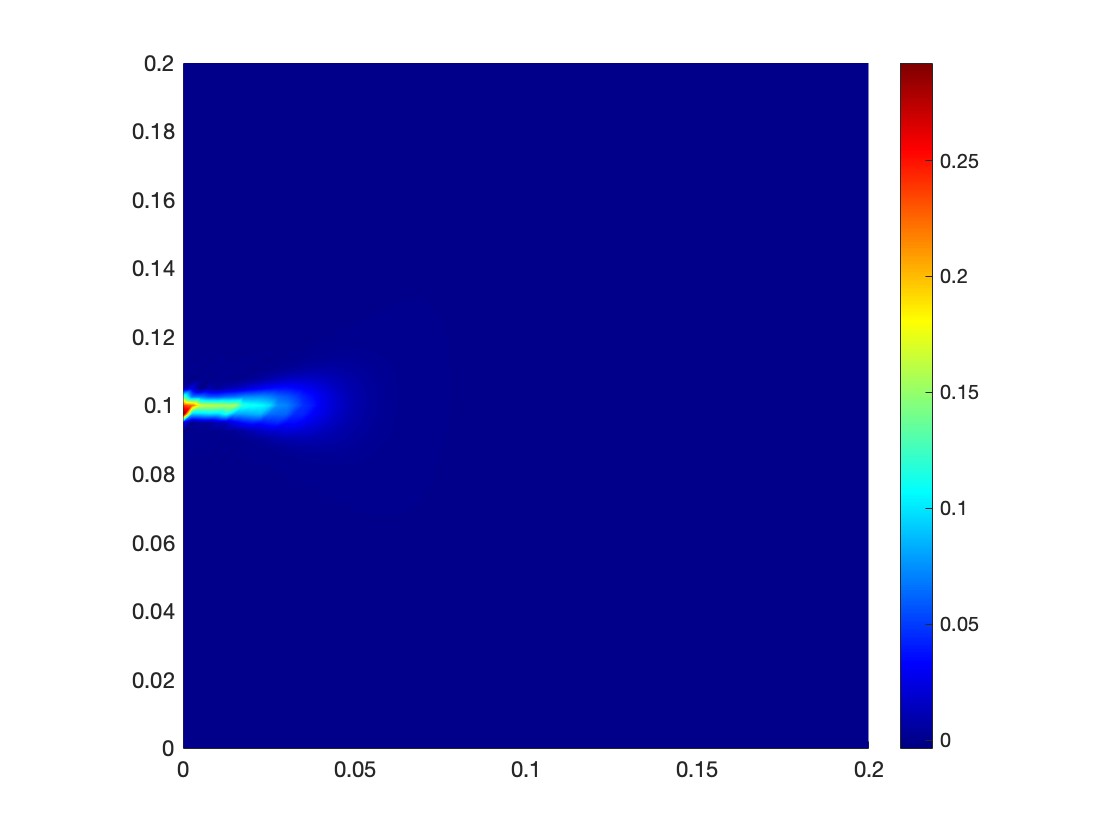}}
\subfigure[]{
    \label{fig2:c}
\includegraphics[width=6.5cm,height=5cm]{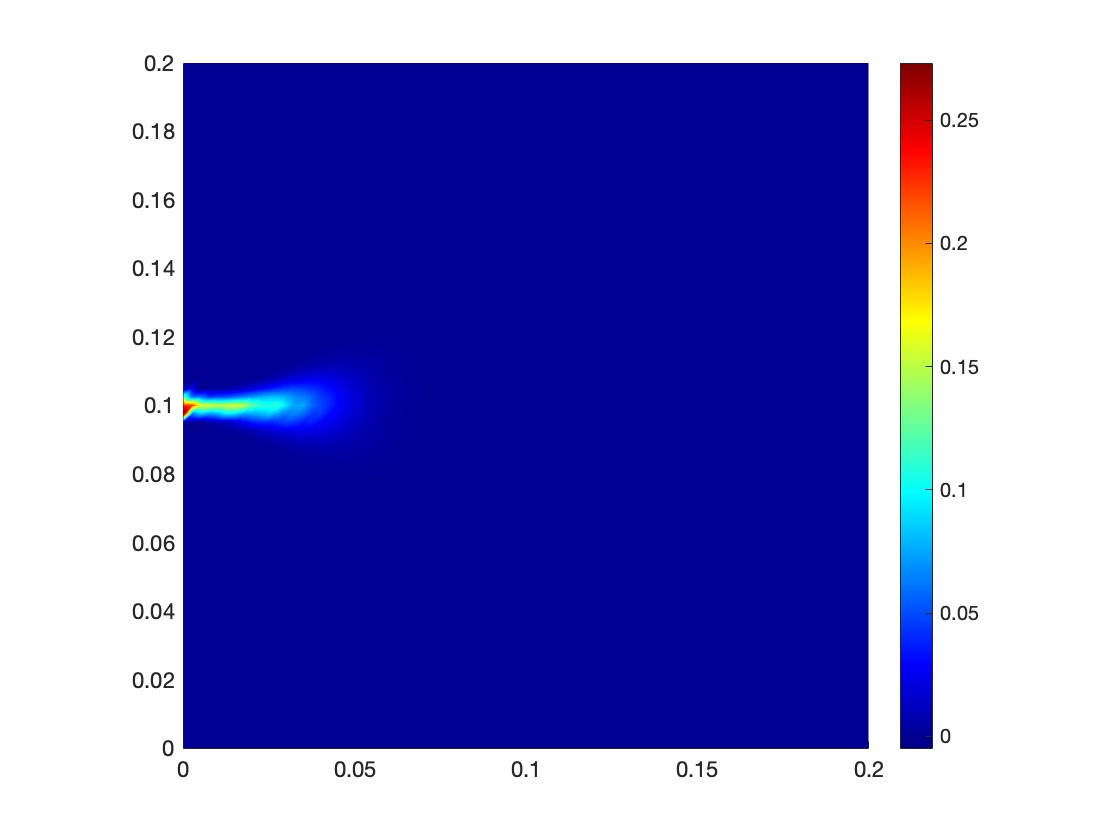}}
\subfigure[]{
    \label{fig2:d}
\includegraphics[width=6.5cm,height=5cm]{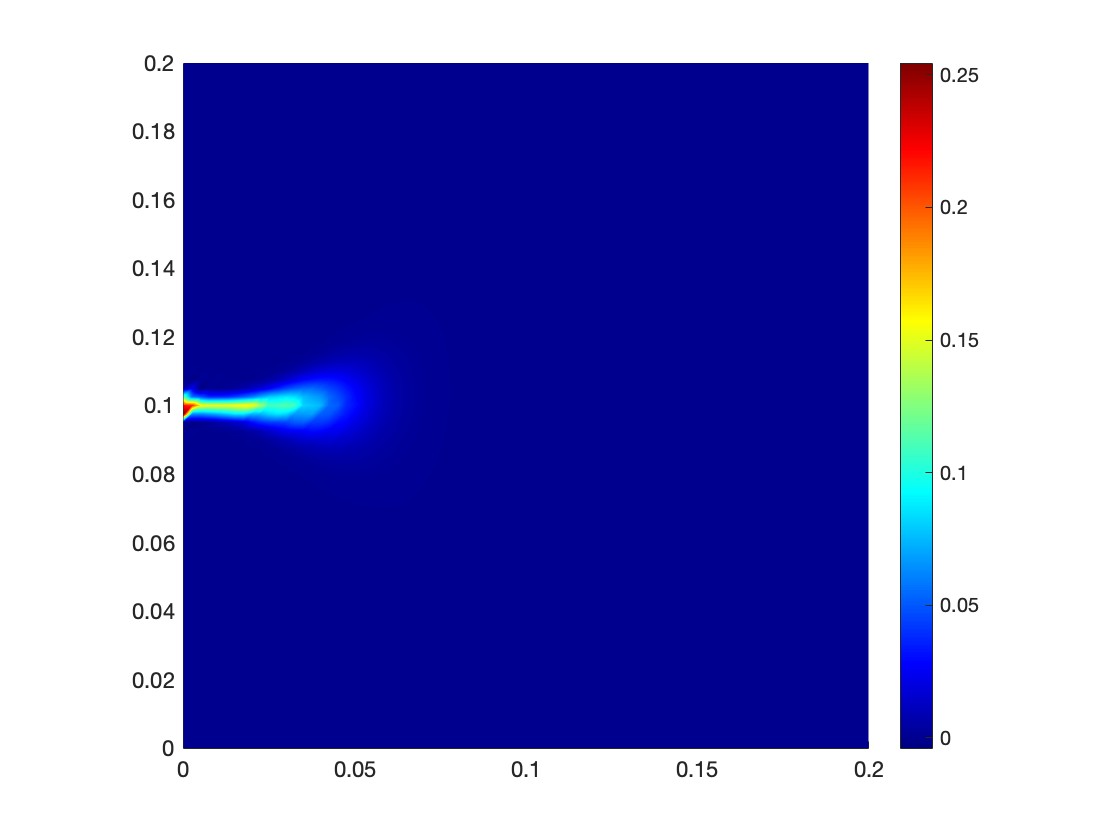}}
\caption{Concentration of acid at time $t=10s$ (a), $t=20s$ (b), $t=30s$ (c),
$t=40s$ (d) with $d=1e-5$.}
\label{fig2}
\end{figure}

The distributions of acid concentration and rock porosity at the final times $t=10s$, $20s$, $30s$ and $40s$ are computed using the parameters $h=1/80$ and $\Delta t=0.1$. These results are illustrated in Figures \ref{fig1}–\ref{fig2}. From the figures, the phenomenon of wormhole propagation is clearly observed.

\bigskip

\section*{Acknowledgments}

 J. Zhang and H. Guo's work was supported partially by the National Key Research and Development Program of China (Grant No. 2023YF1009003), the Natural Science Foundation of Shandong Province (Grant No. ZR2023MA081) and the Fundamental Research Funds for the Central Universities (Grant No. 22CX03020A). J. Zhu’s work was partially supported by the National Council for Scientific and Technological Development of Brazil (CNPq).


\end{document}